\documentclass[reqno,12pt]{amsart}

\usepackage[toc]{appendix}
\usepackage[utf8]{inputenc}
\usepackage[T1]{fontenc}
\usepackage{microtype} 
\usepackage{cite}
\usepackage{graphicx}
\usepackage[caption=false]{subfig}
\usepackage{hyperref}
\usepackage{amsfonts,amsmath,amsthm,amssymb}
\usepackage[a4paper,bindingoffset=0.5cm,left=2cm,right=2cm,top=2.5cm,bottom=2cm,footskip=.8cm]{geometry}

\usepackage{tikz}
\usetikzlibrary{automata,positioning}
\usetikzlibrary{decorations.pathreplacing,decorations.markings,shapes.geometric}
\usetikzlibrary{shapes,arrows}
\usetikzlibrary{backgrounds,calc,positioning}

\def \E {{\mathbb E}}
\def \ee {{\rm e}}

\theoremstyle{definition}

\newtheorem{example}{Example}[section]
\theoremstyle{plain}
\newtheorem{theorem}{Theorem}[section]
\newtheorem{proposition}{Proposition}[section]

\newtheorem{remark}{Remark}[section]
\newtheorem{corollary}{Corollary}[section]

\newcommand{\vb}{\vspace{3mm}}

\begin{document}

\allowdisplaybreaks

\title[]{\small 
SINGLE-SERVER QUEUES UNDER OVERDISPERSION\\
IN THE HEAVY-TRAFFIC REGIME}

\author{By O. Boxma, M. Heemskerk, \& M. Mandjes}

\begin{abstract}

This paper addresses the analysis of the queue-length process of single-server queues under overdispersion, i.e., queues fed by an arrival process for which the variance of the number of arrivals in a given time window exceeds the corresponding mean. 
Several variants are considered, using concepts as mixing and Markov modulation, 
resulting in different models with either endogenously triggered or 
exogenously triggered random environments.
Only in special cases explicit expressions can be obtained, e.g.\ when the random arrival and/or service rate can attain just finitely many values. 
While for more general model variants exact analysis is challenging, one {\it can} derive limit theorems in the heavy-traffic regime. 
In some of our derivations we rely on evaluating the relevant Laplace transform in the heavy-traffic scaling using Taylor expansions, whereas other results are obtained by applying the continuous mapping theorem.

\vb

\noindent
{\sc AMS Subject Classification (MSC2010).} 
Primary: 
60K25; 
Secondary: 
60J60, 
90B22 

\vb

\noindent
{\sc Affiliations.} Onno Boxma is with E{\sc urandom} and the Department of Mathematics and Computer Science; Eindhoven University of Technology; P.O. Box 513, 5600 MB Eindhoven; The Netherlands ({\it email}: {\tt \scriptsize o.j.boxma@tue.nl}). Mariska Heemskerk and Michel Mandjes are with Korteweg-de Vries Institute for Mathematics, University of Amsterdam; Science Park 904, 1098 XH Amsterdam; The Netherlands ({\it email}: {\tt\scriptsize  {j.m.a.heemskerk|m.r.h.mandjes}@uva.nl}). Version: \today.

\vb

\noindent
{\sc Acknowledgments.}   The research of all authors is partly funded by the NWO Gravitation Programme N{\sc etworks} (Grant Number 024.002.003). The research of the first and last author is partly funded by an NWO Top Grant (Grant Number 613.001.352). 
 
\end{abstract}

\maketitle

\section{Introduction}

After the publication, fifty years ago, of {\em The Single Server Queue} \cite{Cohen69},
much research effort has been devoted to 
relaxing  modeling assumptions that are commonly imposed in queueing theory. 
Over the past decades this led to various generalizations of the standard 
single-server queue, thus getting a handle on more versatile and realistic model variants.

Arguably the most prominent generalization concerned the introduction of {\it Markov modulation} \cite{Neuts}.
Under this paradigm, one can depart from the usual assumption that the arrival and service processes are of the renewal type.  
In its most elementary form, the distributions of the interarrival times and service times 
change at transition epochs of a Markovian background process; in addition, there are variants where this happens at arrival or service completion epochs.
A key model is the Markov-modulated M/M/1 queue: when the background process moves to, say, state $j$, the arrival rate becomes $\lambda_j$ while the service rate becomes $\mu_j$. 
Service modulation comes in different flavors: (i) {\it customer-based}, (ii)~{\it server-based}, and what could be called (iii)~{\it real-time updated} service modulation. 
In the first case the state of the background process determines the service requirement of a customer upon arrival, whereas in the second case the service duration is determined by the background state at service initiation.
However, only models with real-time updated service rates fit in the so-called {\it Quasi-Birth-Death} (QBD) framework: they are described by Markov chains defined on a level-phase state space, where level transitions are skipfree.

Queues with a Quasi-Birth-Death structure form a particularly interesting class of Markov-modulated queues.
Triggered by the pioneering work of Neuts \cite{Neuts}, and thanks to major contributions by e.g.\ Latouche, Ramaswami, and Taylor,  
the theory of QBD queues is well developed. 
It provides elegant matrix expressions for key performance measures which can be evaluated through efficient numerical algorithms. 
For general overviews we refer to the book \cite{LatoucheRamaswami} and the proceedings \cite{LT1,LT2}.
An important generalization is provided by {\it queues with Markov additive input}, which can be viewed as Markov-modulated L\'evy processes that are reflected at 0;
we refer to  \cite[Ch.\ XI]{Asmussen} for an authoritative exposition.

\vb

Recently the single-server queue has been generalized in other directions as well.
Various data studies  indicated  that often the variance of the number of arrivals in a given time window exceeds the corresponding mean, a phenomenon often referred to as {\it 
overdispersion}.  
This observation led one to question the standard assumption of Poisson arrivals (under which the variance would equal the mean). 
In order to create overdispersion in a model for call center traffic,
in \cite{JongbloedKoole} the authors suggest to use a Poisson {\em mixture model} for the number of arrivals in an isolated time slot:
the deterministic arrival rate $\lambda$ of the Poisson distribution is replaced by a random variable $\Lambda$.

Following up on the idea of a random arrival rate, the standard Poisson arrival process as a whole can be replaced by a Coxian arrival process, i.e., by a Poisson process for which the intensity  itself is a stochastic process $\Lambda(\cdot)\equiv \{\Lambda(t), t \geq 0\}$.
Infinite-server queues with Coxian input can be analyzed in great detail, essentially due to the property that individual customers do not interfere with each other.
The most elementary variant is the one in which the intensity is resampled, in an i.i.d.\ manner, at equidistant epochs \cite{HLM17}. 
This framework is extended in \cite{HLMM20}, allowing the sample rates to be dependent in an autoregressive manner. 
In \cite{Koops} the object of study is the infinite-server queue with $\Lambda(\cdot)$ corresponding to a `shot-noise' intensity process, whereas in other recent papers $\Lambda(\cdot)$ is a (generalized) Hawkes, or `self-exciting' process \cite{DawPender,DP2, Koops2}. 
While infinite-server queues with Coxian input allow explicit analysis, this usually does not hold for their single-server counterparts. 
See \cite{AlbrecherAsmussen} for an early contribution to the analysis of a specific
insurance risk model (which can be considered dual to a single-server queue) with a Coxian claim arrival process; the authors derive large-deviations results for the ruin probabilities.

The aim of the present paper is to explore single-server queues under overdispersion, where this overdispersion is realized by imposing specific mixing mechanisms. 
It turns out that for these models exact analysis is in general highly challenging. 
This is why we resort to the heavy-traffic scaling, in which the system load approaches unity. 
In this regime for various model variants explicit  limiting results {\it can} be derived.

\vb

\noindent
{\it Related literature -- Markov modulation and mixing}.
Above we already gave a few key references on QBD processes and Markov-modulated queues.
In addition, we would like to mention a paper of Regterschot and de Smit \cite{RS},
making an important methodological contribution by developing a matrix Wiener-Hopf approach to study the M/G/1 queue with Markov-modulated arrivals and service requirements (where the service modulation is customer-based). 
Another notable contribution is by Prabhu and Zhu \cite{PZ}, who work with a very similar but more comprehensive model, using techniques based on infinitesimal generators to analyze the waiting time, idle time, and busy period.

An early example of mixing in the queueing literature is provided in \cite{Abate1}.
Here a GI/G/1 queue is considered, with the service times being exponentially distributed, but with a mean that is a Pareto distributed random variable.
The attractive feature of this construction is that the resulting service-time distribution has an explicit Laplace transform,
even though the distribution is heavy-tailed. A similar procedure has been followed in \cite{Cohen97}:
for Pareto distributed service times with one of its parameters being Gamma distributed,  a mathematically convenient Laplace transform was identified. 
This class of distributions was further generalized in \cite{AW}, considering two classes of so-called Beta mixtures of exponentials.
In \cite{RAB}, the waiting-time distribution of a single-server queue is analyzed for models in which the arrival rate, the service rate
or the traffic load, is random. The paper also covers a duality result between such queueing models and a class of insurance risk models, thus allowing one (i)~to obtain some new results for insurance risk models in which a parameter is random,
and (ii)~to translate some insurance risk results, with mixing, from \cite{Constantinescu} to queueing results. 
Other references on single-server queues with random arrival and/or service rates are \cite{FSW,FR}.

In the references featuring in the previous paragraph, the random parameter was sampled once and for all. In 
\cite{Klaasse} other sampling procedures are explored; specifically, 
if parameters are resampled after each busy period, then quite detailed results can be obtained.

\vb

\noindent
{\it Related literature -- heavy-traffic}.
Heavy-traffic analysis has a tradition within the queueing literature that goes back to the 1960s.
It concerns the study of scaled random quantities within the framework of queueing systems, in the asymptotic regime where the queue's load is increasingly heavy.
One of its pioneers, Kingman \cite{King1,King2}, derived the asymptotic distribution of the scaled steady-state sojourn time and  queue length (and the waiting time) in a GI/G/1 queue; both quantities converge to exponentially distributed random variables, under the condition that the interarrival times and service times have finite variance.
The proof is classical: in the corresponding transform the load of the system is increased to 1, so as to obtain the transform of the exponential distribution in the limit. 
Later also path-level heavy-traffic limit results were established: when scaling time as well, usually relying on the continuous mapping theorem, it was shown that the queue-length process weakly converges to reflected Brownian motion. 
Such limit theorems are particularly useful for queues that do not allow an explicit performance analysis, bearing in mind that, conveniently, for the limiting (Brownian-motion related) objects a broad range of closed-form expressions is known. 
In \cite{WHI} an overview, covering various heavy-traffic limit results, is given, including an account of its relation to the vast literature on diffusion approximations. 
We also refer to the textbook treatments in \cite[Section X.7]{Asmussen} and \cite[Ch.\ V]{DM}. 
Heavy-traffic analyses for Markov-modulated single-server queues {and their QBD counterparts} are given in e.g.\ \cite{ASMHT,BurmanSmith, DIM, TV}. 

\vb

\noindent
{\it Main results and organization of the paper}.
Section~\ref{sec:exact} sketches an approach to analyze the M/M/1 queue in which the arrival rate $\Lambda$ and the service rate ${\mathcal M}$ are resampled at i.i.d.\ exponentially distributed intervals.
The main goal of the section is to point out that exact analysis is within reach when $\Lambda$ and  ${\mathcal M}$ can attain only finitely many values, but complications arise otherwise. 
With these complications in mind, in the remainder of the paper we mainly focus on heavy-traffic analysis. 
In the next section we analyze three different overdispersed queueing models, each with its own resampling mechanism. 

In Section~\ref{sec:exo} we consider the most basic model: an M/M/1 queue where the rate vector $(\Lambda,{\mathcal M})$ is resampled (in an i.i.d.\ manner) at Poisson epochs, with $(\Lambda,{\mathcal M})$ attaining only finitely many values. It takes a little thought to realize that this model is a special case of the Markov-modulated M/M/1 queue, with the Markovian background process having a finite state space.
For this more general model we first show how the Laplace transform of the steady-state distribution of the number of customers $Q$, jointly with the background state $J$, can be determined (Section~\ref{sec:3.1}).
Then in Section~\ref{HTSL} we study the heavy-traffic scaling limit of $Q$: by letting the traffic load $\rho$ tend to $1$ in the Laplace transform we prove that $(1-\rho)Q$ converges to an exponentially distributed random variable.
In Section~\ref{sec:3.3} we narrow our scope to the case where $(\Lambda,{\mathcal M})$ are resampled, under which the parameter of the limiting exponential distribution simplifies considerably. 

Section~\ref{sec:exoperiodic} is also concerned with an M/M/1 queue with resampling
of the arrival rate and service rate, but now the support of $(\Lambda,{\mathcal M})$ is not restricted to finitely many values, and the resampling intervals are not necessarily exponentially distributed.
We first determine the mean, variance and covariance of the cumulative arrival process and cumulative potential service process.
Those results are subsequently used to show, using the continuous mapping theorem, that an appropriately scaled version of the queue-length process converges to reflected Brownian motion
when the traffic load approaches unity. {Our result does {\it not} directly imply that the stationary number of customers under the heavy-traffic scaling converges to the stationary version of reflected Brownian motion (which has an exponential distribution).}

While the resampling mechanisms featuring in Sections~\ref{sec:exo} and \ref{sec:exoperiodic} can be seen as exogenously triggered, the resampling in Section~\ref{sec:endo} is endogenous: we consider an M/G/1 queue in which at every service completion the arrival rate is resampled.
After providing an exact analysis of the transient queue-length distribution right after service completions, we use the obtained results to show that a scaled version of the {transient} queue length (where the scaling involves both space and time) converges {to its counterpart for reflected Brownian motion} when the traffic load $\rho$ approaches $1$.
Also in this context we succeed in proving that the {scaled stationary queue length} $(1-\rho)Q$ converges to an exponentially distributed random variable.

We conclude the paper in Section~\ref{sec:conclusion} with a brief discussion and some suggestions for further research. 
Importantly, the Sections \ref{sec:exo}--\ref{sec:endo} focus on three different models; for this reason, we introduce for each of them specific notation in the corresponding section.

\section{Exact results: an exploration} \label{sec:exact}
In this section we consider an M/M/1 queue with the special feature that at Poisson epochs, the underlying  rate vector $(\Lambda,{\mathcal M})$, which is componentwise non-negative, is resampled from a given distribution. 
In particular, we do not assume that $\Lambda$ and ${\mathcal M}$ are independent.
Let $q^{-1}$ be the value of the mean `inter-sample time'. 
The main objective of this section is to examine to what extent the resulting queueing system allows a closed-form solution.

Starting point are the following results from \cite[Section I.4.4]{Cohen69}
for the transient behavior of a birth-death process with constant birth and death rates.
Let $Q_t$ denote the number of customers in the M/M/1 queue at time $t \geqslant 0$.
Assume that, at time $0$ and until $\xi \sim$ exp($q$), the rates are given by $\Lambda_0 = \lambda$ and 
${\mathcal M}_0 = \mu$, with traffic load $\rho := \lambda/\mu$.
To describe the number of customers at the random epoch $\xi$, the following quantities play a role \cite[Section I.4.4]{Cohen69}.
Let $x_1\equiv x_1(\lambda,\mu,q)$ and $x_2\equiv x_2(\lambda,\mu,q)$ be the roots of 
$F(z)\equiv F(z\,|\,\lambda,\mu,q):=\lambda z^2- (\lambda + \mu + q)z+ \mu$, where $x_1$ is the larger one:
\begin{equation}
x_{1} = \frac{\lambda + \mu + q + \sqrt{(\lambda + \mu + q)^2 - 4 \lambda \mu }}{2 \lambda},\:\:\:\:
x_{2} = \frac{\lambda + \mu + q - \sqrt{(\lambda + \mu + q)^2 - 4 \lambda \mu }}{2 \lambda} .
\end{equation}
Then, according to \cite[(4.27) on p.\ 80]{Cohen69} (taking into account the scaling with respect to $\mu$ mentioned there),
for $|z| \leqslant 1$ and $i=1,2,\ldots$,
\begin{eqnarray} 
{\mathbb E}\big(z^{Q_\xi}|(Q_0,\Lambda_0,\mathcal{M}_0)=(i,\lambda, \mu)\big)
&=& \sum_{j=0}^{\infty} z^j \int_{t=0}^{\infty} q {\rm e}^{-q t} {\mathbb P}(Q_t=j|(Q_0,\Lambda_0,\mathcal{M}_0)=(i,\lambda, \mu)) {\rm d}t
\nonumber
\\
&=& \frac{q}{\mu} \frac{(1-z)x_2^{i+1} - (1-x_2)z^{i+1}}{(1-x_2)(\rho z^2 -(1 + \rho +  {q}/{\mu})z+1)} 
\nonumber
\\
&=& q \frac{(1-z)x_2^{i+1} - (1-x_2)z^{i+1}}{(1-x_2)F(z)} .
\label{Cohen1}
\end{eqnarray}
As is well-known \cite[p.\ 190]{Cohen69},
the smaller root $x_2$ can be interpreted as the Laplace-Stieltjes transform ${\mathbb E}[{\rm e}^{- q P}]$ of the busy period $P$ in an M/M/1 queue
with arrival rate $\lambda$ and service rate $\mu$, given that $\lambda < \mu$.

From now on we assume that $\E[\Lambda] < \E[\mathcal M]$, implying that the system is stable.
If the system was already in steady state at time $0$, and with $Q$ denoting the steady-state queue length,
then we can write $Q_0 = Q$ and
\begin{equation}
{\mathbb E}[z^{Q_\xi} \mid (\Lambda_0, {\mathcal M}_0) = (\lambda, \mu)]
=  \frac{q (1-z)x_2 \,\E[x_2^Q]}
{(1-x_2)F(z)} 
-  \frac{q z\, \E[z^Q]}
{F(z)}.
\label{Cohen2}
\end{equation}
Note that the queue length $Q_\xi$ at the resample epoch $\xi$ has the same steady-state distribution as $Q_0$, as at time $\xi$ rates $(\Lambda_1, {\mathcal M}_1)$ are resampled in an i.i.d.\ fashion from the same distribution as $(\Lambda_0,{\mathcal M}_0)$.
Hence we can obtain an expression for the probability generating function (PGF) of $Q$ from (\ref{Cohen2}) by integrating both sides with respect to $\lambda$ and $\mu$ (recalling the dependence of $F(z)$ and $x_2$ on $\lambda$ and $\mu$):
\begin{align}
{\mathbb E}[z^{Q}]
&={q (1-z)} \int_{\lambda = 0}^{\infty} \int_{\mu =0}^{\infty} \frac {x_2\E[x_2^Q]} {(1-x_2) F(z)}{\mathbb P}(\Lambda \in {\rm d} \lambda,{\mathcal M} \in {\rm d} \mu)
\nonumber
\\
&\hspace{4mm}- q z\, \E[z^Q]
\int_{\lambda = 0}^{\infty} \int_{\mu =0}^{\infty} 
\frac{1}{F(z)}  
{\mathbb P}(\Lambda \in {\rm d} \lambda,{\mathcal M} \in {\rm d} \mu) .
\label{Cohen3}
\end{align}
Hence,
\begin{eqnarray}
\E[z^{Q}]
&=& q (1-z) \left(\frac{\displaystyle \int_{\lambda = 0}^{\infty} \int_{\mu =0}^{\infty}  \frac{x_2 \E[x_2^Q]}
{\displaystyle(1-x_2)F(z)} {\mathbb P}(\Lambda \in {\rm d} \lambda,{\mathcal M} \in {\rm d} \mu)}{\displaystyle
1+ q z
\int_{\lambda = 0}^{\infty} \int_{\mu =0}^{\infty} 
\frac{1}{F(z)}   
{\mathbb P}(\Lambda \in {\rm d} \lambda,{\mathcal M} \in {\rm d} \mu)} \right).
\label{Cohen4}
\end{eqnarray}
The crucial insight is that a major complication arises from the fact that the numerator in the right-hand side of (\ref{Cohen4}) still contains an unknown expression.
More than that, since $x_2$ is a function of $\lambda$ and $\mu$,
there may be multiple
unknowns $\E[x_2^Q]$;
even infinitely many as $\Lambda$ and $\mathcal{M}$ need not live on a finite state space.

To study this complication in greater detail, we proceed by briefly pointing out the crucial difference between the case in which
the vector $(\Lambda,{\mathcal M})$ has a finite support (i.e., can attain finitely many values),
and the case in which this vector can take on infinitely many values.
First consider the case that $(\Lambda,{\mathcal M})$ can attain two values, i.e., we assume that, for $\pi_1\in(0,1)$,
\begin{equation}
{\mathbb P}(\Lambda = \lambda_1,{\mathcal M}=\mu_1) = \pi_1, ~~~
{\mathbb P}(\Lambda = \lambda_2,{\mathcal M}=\mu_2) = \pi_2 = 1 - \pi_1.
\label{Cohen5}
\end{equation}
In this case,
writing $x_{2i}$ to indicate its dependence on the values of $\lambda_i$ and $\mu_i$, and with
$F_i(z):=\lambda_i z^2- (\lambda_i + \mu_i + q)z+ \mu_i$,
(\ref{Cohen4}) reduces to
\begin{eqnarray}
\E[z^{Q}]
&=&
q(1-z)\left(\frac{\displaystyle \sum_{i=1}^2  \frac{\pi_i\,x_{2i} \E[x_{2i}^Q]}
{(1-x_{2i})F_i(z)}}
{\displaystyle
1+ q z
\sum_{i=1}^2
\frac{\pi_i}{F_i(z)}}\right),
\label{Cohen6}
\end{eqnarray}
the only unknowns being $\E[x_{2i}^Q]$ for $i=1,2$.
After multiplication by
$F_1(z)\,F_2(z)$,
the denominator in the right-hand side of (\ref{Cohen6}) is a polynomial $D(z)$ in $z$ of order $4$,
with a zero at $z=1$ and three other zeros which can be explicitly obtained via Cardano's formula.
When the stability condition
$\pi_1\lambda_1 + \pi_2 \lambda_2 < \pi_1 \mu_1 + \pi_2 \mu_2$ holds,
one of those zeros lies in $(0,1)$. 
In fact, it is easily verified that $D(0) = \mu_1 \mu_2 > 0$, $D(1)=0$, and $D'(1) >0$ iff that stability condition holds.
The  zero in $(0,1)$ and the zero $z=1$ should also be zeros of the numerator in (\ref{Cohen6}),
which gives rise to two linear equations in the two above-mentioned unknowns.
We conclude that this problem is solvable.

When the vector $(\Lambda,{\mathcal M})$ can attain a larger, but still finite, number of values,
one typically aims for a Rouch\'e-type argument to prove that there is a specific number of zeros in $|z| \leqslant 1$.
We do not discuss this further, because
this model can be viewed as a special case of a Markov-modulated M/M/1 queue,
which has been analysed in detail, cf.\ \cite{RS}.
Via an intricate argument, it is there proven for the case of a background process with $d$ states
(and while studying the waiting time rather than queue length)
that a particular determinant has, in steady state, exactly $d-1$ zeros in the right-half plane and one at zero.
That knowledge can be exploited to determine the vector of waiting-time transforms
(of length $d$, distinguishing between waiting times of customers arriving in different background states).

However, this approach breaks down when the state space of the background process is not finite.
In particular, when the support of $(\Lambda,{\mathcal M})$ is uncountable, the numerator in the right-hand side of (\ref{Cohen4})
is an unknown function of $z$, and the denominator may have zeros on a contour in $|z| \leqslant 1$.
It is far from obvious how to exploit that knowledge to find the numerator in (\ref{Cohen4}).

The main conclusion of this section is that, in order to obtain exact results, complications arise when the support of $(\Lambda,{\mathcal M})$ is not finite. Classical queueing-theoretic techniques fail, and new approaches need to be developed. In the next sections we focus on the heavy-traffic regime, in which rather explicit results {\it can} be obtained.

\section{Exogenously triggered resampling, finite support}
\label{sec:exo}
In this section we consider a heavy-traffic analysis of a Markov-modulated M/M/1 queue with finitely many background states. 
It covers the setting discussed in Section $\ref{sec:exact}$: an M/M/1 queue in which at Poisson epochs the arrival rate and service rate are resampled from a distribution with finite state space. 
As turns out in Section $\ref{sec:3.3}$, in this special case the (parameter of the) heavy-traffic limiting distribution simplifies considerably.

\subsection{Model description}
\label{sec:3.1}
In this section we consider the following Markov-modulated M/M/1 queue. 
There is a background process, which is assumed to be irreducible, and which jumps from $i$ to $j\not=i$ with rate $q_{ij}$, with $i,j\in\{1,\ldots,d\}$. 
We define $q_i:=-q_{ii}:=\sum_{j\not=i}q_{ij}$. 
When the background state is $i$, the arrival rate is $\lambda_i\geqslant 0$ and the service rate is $\mu_i\geqslant 0$. 
As before, $Q_t$ denotes the queue length at time $t$ and $Q$ the steady-state queue length.

As a first step, we characterize the PGF of $Q$. 
Let $p_t(k,i)$ be the probability that at time $t$ there are $k$ customers (with $k\in{\mathbb N}_0$) and the background state is $i$ (with $i\in\{1,\ldots,d\}$).
As $\Delta\downarrow 0$, we obtain by a classical Markovian argumentation, for any $k\in{\mathbb N}_0$,
\begin{align*}
p_{t+\Delta}(k,i)&=\lambda_i\Delta \,p_t(k-1,i)1_{\{k\geqslant 1\}} + \mu_i\Delta\,p_t(k+1,i) +\sum_{j\not=i}q_{ji}\Delta\,p_t(k,j)\:+\\
&\:\:\:\:\:\:\:\hspace{9mm}\big(1-(\lambda_i+\mu_i1_{\{k\geqslant 1\}}+q_i)\Delta\big)\,p_t(k,i)+o(\Delta).
\end{align*}
Multiplying by $z^k$ and summing over $k$ yields, in a standard manner, a relation in terms of the PGF of $Q_{t+\Delta}$ and $Q_t$, jointly with the corresponding states of the background process. 
By subtracting ${\mathbb E}[z^{Q_{t}}1_{\{J_{t}=i\}}]$ from both sides of this relation, dividing by $\Delta$, and letting $\Delta\downarrow 0$, we obtain the differential equation
\begin{align*}\frac{{\rm d}}{{\rm d}t}{\mathbb E}[z^{Q_{t}}1_{\{J_{t}=i\}}] &= 
\lambda_i(z-1)\ {\mathbb E}[z^{Q_{t}}1_{\{J_{t}=i\}}] \:+\\
&\:\:\:\:\:\:\:\mu_i\left(\frac{1}{z}-1\right){\mathbb E}[z^{Q_{t}}1_{\{J_{t}=i, Q_t>0\}}]+
\sum_{j=1}^d q_{ji}\ {\mathbb E}[z^{Q_{t}}1_{\{J_{t}=j\}}].
\end{align*}
We send $t$ to $\infty$ to obtain a system of differential equations for the steady-state queue length $Q$, jointly with the state $J$ of the background process. 
Denoting the corresponding PGF by $f_i(z):= {\mathbb E}[z^{Q}1_{\{J=i\}}]$,
it thus follows that
\begin{align}\label{betaeq}
0&= 
\lambda_i(z-1)\,f_i(z) +
\mu_i\left(\frac{1}{z}-1\right)\big(f_i(z)-\beta_i\big)+
\sum_{j=1}^d q_{ji}\,f_j(z),
\end{align}
where $\beta_i$ denotes ${\mathbb P}(Q=0,J=i)$. 
This equation can, for any $z\in(0,1)$ and for a given vector of probabilities ${\boldsymbol\beta}$, be rewritten as a system of linear equations. 
Concretely, with the $(i,j)$-th entry of the matrix $A(z)$ given by 
\[a_{ij}(z):=\lambda_i\,(z-1)1_{\{i=j\}} +\mu_i\left(\frac{1}{z}-1\right)1_{\{i=j\}}+q_{ji},\] 
and the $i$-th entry of the vector 
${\boldsymbol b}(z\,|\,{\boldsymbol \beta})$ given by 
\[b_i(z\,|\,{\boldsymbol \beta}):=\mu_i\left(\frac{1}{z}-1\right)\,\beta_i,\] 
we arrive at the equation
$A(z)\,{\boldsymbol f}(z) = {\boldsymbol b}(z\,|\,{\boldsymbol \beta}).$
We thus obtain, modulo the invertibility of the matrix $A(z)$,
\[{\boldsymbol f}(z) =(A(z))^{-1}\, {\boldsymbol b}(z\,|\,{\boldsymbol \beta}).\]
By Cramer's rule, we have that, with $A_i(z)$ defined as $A(z)$ but with the $i$-th row replaced by ${\boldsymbol b}(z\,|\,{\boldsymbol \beta})$, and with
$\alpha_i(z):={\rm det}\,A_i(z)$ and $\alpha(z):={\rm det}\,A(z)$,
\begin{equation} \label{fieq}
f_i(z) =\frac{\alpha_i(z)}{\alpha(z)}.
\end{equation}
We have thus found the PGF of $Q$, jointly with the background state $J$, in terms of the vector of probabilities ${\boldsymbol \beta}$. 
Using the fact that zeros of the denominator of (\ref{fieq}) should be zeros of the corresponding numerator as well, these probabilities can in principle be found. As it will turn out below, however, their precise value does not affect the heavy-traffic results.

\subsection{Heavy-traffic scaling limit}\label{HTSL}
Now that we have an expression for the PGF of $Q$, this can be exploited to get insight into its heavy-traffic behavior. 
Following the ideas of \cite[Ch.\ VI]{ABI}, we do so by expanding ${\boldsymbol f}(z)$ in $z=\ee^{-(1-\rho)s}$ as $\rho\uparrow 1$, with $\rho := {\boldsymbol \pi}{\boldsymbol \lambda}/ {\boldsymbol \pi}{\boldsymbol \mu}$.
Here ${\boldsymbol \pi}$ denotes the row vector of steady-state probabilities for the background process.

\begin{remark}\label{rem:scaling}{\em 
We let the heavy-traffic regime be reached by scaling ${\boldsymbol\lambda}$ in a uniform fashion, in the sense that a vector of fixed arrival rates ${\boldsymbol \lambda_0}$ is multiplied by a constant, namely by $c>0$ such that $\rho = {\boldsymbol \pi}{\boldsymbol \lambda}/ {\boldsymbol \pi}{\boldsymbol \mu}= {\boldsymbol \pi}(c{\boldsymbol \lambda_0})/ {\boldsymbol \pi}{\boldsymbol \mu}= c \rho_0 \uparrow 1$ as $c \uparrow \rho_0^{-1}$, with $\rho_0 :=  {\boldsymbol \pi}{\boldsymbol \lambda}_0/ {\boldsymbol \pi}{\boldsymbol \mu} < 1$.
Note that the results derived in this section are not affected by this choice, as they only depend on ${\boldsymbol \mu}$, which is kept fixed.
A similar procedure has been followed in Sections \ref{sec:3.3}, \ref{sec:WCBM} and \ref{sec:endoHT}.
}
\end{remark}

The underlying idea is that we show that, as $\rho\uparrow 1$, ${\mathbb E}[\ee^{-s(1-\rho)Q}]$ converges to the Laplace-Stieltjes transform of an exponentially distributed random variable (with a specific rate). 
By L\'evy's convergence theorem this thus implies that $(1-\rho)Q$ converges in distribution to this exponentially distributed random variable.

As a first step, we study the behavior of $\alpha(z)$ and $\alpha_i(z)$ as $z\uparrow 1.$ 
\begin{itemize}
\item[$\circ$]
Note that $a_{ij}(1)=q_{ji}$; due to the singularity of the transition rate matrix $(q_{ij})_{i,j=1}^d$, it is directly seen that $\alpha(1)=0$. 
\item[$\circ$]
As $f_i(z)\in[0,1]$ for all $z \in (0,1)$, we conclude that also $\alpha_i(1)=0$, for $i=1,\ldots,d$. 
\item[$\circ$]
To study $\alpha'(1)$, we first observe that $\alpha(z)={\rm det}\,\bar A(z)$, with $\bar A(z)$ defined as $A(z)$ with the first row replaced by the sum of all rows; i.e., the $(1,j)$-th entry is $\gamma_j(z):=\lambda_j(z-1)+\mu_j(1/z-1)$ (recalling that the rows of the transition rate matrix sum to $0$). Using standard rules for evaluating determinants, \[\alpha(z)=\sum_{j=1}^d \gamma_j(z) \,\delta_j(z),\] where $\delta_j(z)$ is the determinant of the appropriate $(d-1)\times(d-1)$ cofactor matrix. We hence find
\[\alpha'(1) =\sum_{j=1}^d \big(\gamma'_j(1)\,\delta_j(1) +\gamma_j(1)\,\delta'_j(1)\big).\] To evaluate the right-hand side of this equation, we first note that $\gamma_j(1)=0$. Also, $\gamma'_j(1) = \lambda_j-\mu_j.$ We are thus left with evaluating $\delta_j(1)$.

\noindent
We proceed by showing that $\delta_j(1)$ is proportional to $\pi_j$, as follows. 
The vector ${\boldsymbol \pi}$ can be found by solving the linear system of equations $\sum_{j=1}^d \pi_j q_{ji}= 0$ for $i = 1, \dots, d$ (of which one equation is redundant) together with ${\boldsymbol \pi}{\boldsymbol 1} =1$. 
It follows that this system of equations can be written as $T {\boldsymbol \pi}^\top = {\boldsymbol e}_1$, with the $(i,j)$-th element of $T$ being defined by~$1$ if ${i=1}$ and by $q_{ji}$ else, and ${\boldsymbol e}_1$ defining the first unit vector. 
The vector ${\boldsymbol \pi}$ can again be evaluated using Cramer's rule. 
More concretely, with $T_j$ being equal to the matrix $T$ but with the $j$-th column replaced by ${\boldsymbol e}_j$ and $\tau:= {\rm det}\,T$, we find
\[\pi_j = \tau^{-1}\cdot {\rm det}\,T_j.\]
It is immediately verified that ${\rm det}\,T_j=\delta_j(1)=\tau\pi_j$ (recalling that $a_{ij}(1)=q_{ji}$).

\noindent Upon combining the above, we get $\alpha'(1) = \tau\,({\boldsymbol \pi} {\boldsymbol \lambda}- {\boldsymbol \pi} {\boldsymbol \mu})=-\tau\,{\boldsymbol \pi} {\boldsymbol \mu}\,(1-\rho)$. 
\item[$\circ$]
In addition, by L'H\^opital's rule,
\[1=\sum_{i=1}^d {\mathbb P}(J=i)=\lim_{z\uparrow 1} \sum_{i=1}^d {\mathbb E}[z^{Q}1_{\{J=i\}}] = \lim_{z\uparrow 1} \sum_{i=1}^d \frac{\alpha_i(z)}{\alpha(z)}= \sum_{i=1}^d \frac{\alpha'_i(1)}{\alpha'(1)},\]
so that $\sum_{i=1}^d {\alpha'_i(1)} = \alpha'(1) = -\tau\,{\boldsymbol \pi} {\boldsymbol \mu}\,(1-\rho)$ as well. 
Realize that $\tau$ depends on the transition rates $q_{ji}$ only, i.e., not on the arrival rates $\lambda_i$ and service rates $\mu_i$.
\end{itemize}
We thus find that, for a given $s\geqslant 0$, 
\begin{equation*}
\alpha(\ee^{-(1-\rho)s}) = \alpha\big(1 - (1-\rho)s +\tfrac{1}{2}(1-\rho)^2s^2+ O\left((1-\rho)^3\right)\big),
\end{equation*}
which further expands to
\begin{align} 
&\alpha(1)+\alpha'(1)\big( - (1-\rho)s +\tfrac{1}{2}(1-\rho)^2s^2\big) +\tfrac{1}{2}\alpha''(1)(1-\rho)^2s^2+ O\left((1-\rho)^3\right) \nonumber \\ \label{ch7:expand}
&=\tau\,{\boldsymbol \pi}  {\boldsymbol \mu}\,(1-\rho)^2s+\tfrac{1}{2}\alpha''(1)(1-\rho)^2s^2+ O\left((1-\rho)^3\right).
\end{align}
Along the same lines,
\[\sum_{i=1}^d\alpha_i(\ee^{-(1-\rho)s}) = \tau\,{\boldsymbol \pi} {\boldsymbol \mu}\,(1-\rho)^2s+\tfrac{1}{2}\sum_{i=1}^d\alpha_i''(1)(1-\rho)^2s^2+ O\left((1-\rho)^3\right).\]
We conclude that
\begin{equation}
\label{TRA}\lim_{\rho\uparrow 1} {\mathbb E}[\ee^{-s(1-\rho)Q}] =\frac{\tau\,{\boldsymbol \pi}{\boldsymbol \mu}+\tfrac{1}{2} \lim_{\rho \uparrow 1}\sum_{i=1}^d\alpha_i''(1)\,s}{\tau\,{\boldsymbol \pi}{\boldsymbol \mu}+\tfrac{1}{2}\lim_{\rho \uparrow 1}\alpha''(1)\,s}.\end{equation}
As ${\mathbb P}(Q=0)$ vanishes when $\rho\uparrow 1$, we deduce that the expression in (\ref{TRA}) goes to $0$ as $s\to\infty$.
This immediately entails that
$\sum_{i=1}^d\alpha_i''(1)\to 0$ as $\rho\uparrow 1$. 
Recognizing the Laplace transform of the exponential distribution, and defining $\alpha^\circ:=\lim_{\rho \uparrow 1}\alpha''(1)$,
we have proven the following result.  
 
\begin{theorem} \label{ch7:TH1}
Consider the Markov-modulated M/M/1 queue.
As $\rho\uparrow 1$, we have that $(1-\rho)Q$ converges to an exponentially distributed random variable with mean $\alpha^\circ/(2\tau\,{\boldsymbol \pi}{\boldsymbol \mu}).$
\end{theorem}

Upon inspecting the above derivation, we see that we have in addition proven the asymptotic independence of $Q$ and $J$ in the regime where $\rho\uparrow 1$. 
More precisely, as $\rho\uparrow 1$, using that PGFs uniquely characterize their underlying (joint) distribution, 
we have that the bivariate random vector $((1-\rho)Q, J)$ converges to $(\bar Q,\bar J)$, where $\bar Q$ is exponentially distributed with the mean we identified above and $\bar J$ such that ${\mathbb P}(\bar J=i)=\pi_i$, where, remarkably, $\bar Q$ and $\bar J$ are independent; cf.\ the results in e.g. \cite{ASMHT,DIM,TV}.
Hence we have the following result.

\begin{corollary} {\em
As $\rho\uparrow 1$,  $((1-\rho)Q, J)$ converges to $(\bar Q,\bar J)$, where $\bar Q$ is exponentially distributed with mean $\alpha^\circ/(2\tau\,{\boldsymbol \pi}{\boldsymbol \mu})$ and $\bar J$ has distribution ${\mathbb P}(\bar J=i)=\pi_i$, where  $\bar Q$ and $\bar J$ are independent.
}\end{corollary}
 
\subsection{Heavy-traffic scaling limit in the resampling model}
\label{sec:3.3}
We proceed by considering a special case of a Markov-modulated M/M/1 queue, namely the one that corresponds to resampling the arrival rate and service rate at Poisson epochs, as was introduced in Section~$\ref{sec:exact}$. {We do so following similar steps as in the proof of Theorem \ref{ch7:TH1}, but in the resampling context the parameter of the limiting exponential distribution turns out to simplify considerably.}

We follow the construction introduced in \cite[Section 3]{CDMR}. 
Let ${\boldsymbol\pi}$ be some row vector of probabilities summing to $1$.
Take, for a given $q>0$, the transition rate matrix equal to $q \,{\boldsymbol 1}{\boldsymbol \pi} -qI_d,$ with ${\boldsymbol 1}$ an all-ones vector and $I_d$ the $d$-dimensional identity matrix. 
As follows from the reasoning in \cite{CDMR}, we have thus constructed a model in which the arrival and service rate pair $(\Lambda, {\mathcal M})$ is resampled, in an i.i.d.\ manner, after exponentially distributed times (with mean $q^{-1}$). 
At every resampling time, they attain the values $(\lambda_i,\mu_i)$ with probability $\pi_i$. 
Notably, this construction allows for the $\Lambda$ and ${\mathcal M}$ to be dependent, but their support needs to be finite.
We call the resulting model the {\it resampled M/M/1 queue}. 
 
In Section \ref{HTSL} we determined the heavy-traffic scaling limit by expanding the determinants of $A(z)$ and $A_i(z)$ as $z\uparrow 1$. 
In the special resampling setting defined above, however, it turns out that these expansions take an explicit form.
Define $E(s)=A(\ee^{s})$, i.e., 
\[E(s) := \bar E(s)+q \,{\boldsymbol 1}{\boldsymbol \pi},\:\:\:\bar E(s):=
-qI_d + {\rm diag}\{{\boldsymbol \lambda}\}(\ee^s-1) + {\rm diag}\{{\boldsymbol \mu}\}(\ee^{-s}-1) .\]
Due to the fact that $E(s)$ can be written as the sum of a diagonal matrix and a rank-one matrix, its eigenvalues can be characterized as roots of a function from ${\mathbb R}$ to ${\mathbb R}$, as can be seen as follows. 
To this end, we write 
\[{\rm det}(E(s)-\theta I_d) = {\rm det}(\bar E(s)-\theta I_d)\,{\rm det}(I_d+(\bar E(s)-\theta I_d)^{-1} q \,{\boldsymbol 1}{\boldsymbol \pi}).\]
For $A$ and $B$ matrices of dimensions $m\times n$ and $n\times m$, respectively, we have ${\rm det}(I_m-AB)={\rm det}(I_n-BA)$. We thus conclude that
\begin{align*}
{\rm det}(E(s)-\theta I_d) =&\: {\rm det}(\bar E(s)-\theta I_d)\,{\rm det}(1+\,{\boldsymbol \pi}(\bar E(s)-\theta I_d)^{-1} q \,{\boldsymbol 1})\\
=&\: {\rm det}(\bar E(s)-\theta I_d)\left(1-\sum_{i=1}^d \pi_i\frac{q}{q-\lambda_i(\ee^s-1)-\mu_i(\ee^{-s}-1)+\theta}\right).
\end{align*}
We conclude that the eigenvalues  $\theta_1(s)$ up to $\theta_d(s)$ (for a fixed $s$, that is) are the solutions to
\[\frac{1}{q} = \sum_{i=1}^d \pi_i \frac{1}{q-\lambda_i(\ee^s-1)-\mu_i(\ee^{-s}-1)+\theta}=:\Psi_s(\theta).\]
The existence of these (real) eigenvalues follows from the fact that there are poles at $\theta = -q+\lambda_i(\ee^s-1)+\mu_i(\ee^{-s}-1)$, $i=1,\ldots,d$, at which
$\Psi_s(\cdot)$ jumps from  $-\infty$ to $\infty$, whereas $\Psi_s(\cdot)$ converges to $0$ as $\theta\to\pm\infty$. 
This means that $d-1$ of the roots are between two subsequent poles, whereas the largest root is larger than the largest pole. 

\begin{itemize}
\item[$\circ$]
As $s\to 0$, by observing that the locations of all poles converge to $-q$, it follows that all but one eigenvalues, say $\theta_1(s)$ up to $\theta_{d-1}(s)$, converge to $-q$. Also, the largest one, say $\theta_d(s)$, converges to $0$ as $s\to 0.$ 
\item[$\circ$]
We then determine $\theta_j'(1)$ for $j=1,\ldots,d-1$. 
For the moment we assume that the $\kappa_i:=\lambda_i-\mu_i$ are all distinct; the case in which some of the $\kappa_i$ coincide can be dealt with analogously, with slightly more effort.
Without loss of generality we can put the $\kappa_i$ in increasing order, i.e. $\kappa_1 < \dots < \kappa_d$.
Recalling that $\ee^s-1=s+o(s)$ and $\ee^{-s}-1=-s+o(s)$ as $s \rightarrow 0$, it thus follows that, as $s\to0$, 
$\theta_j(s)\in[-q+\kappa_js,-q+\kappa_{j+1}s)$. 
Write $\theta_j(s)=-q +f_js+o(s)$ for $f_j\in[\kappa_j,\kappa_{j+1})$. 
To find the $f_j$, we wish to solve, in the regime where $s\to 0$,
\[\frac{1}{q}= \sum_{i=1}^d \pi_i\frac{1}{(f_j-\kappa_i)s}.\]
Multiplying both sides by $s$, sending $s$ to $0$, and distinguishing between positive and negative terms, we conclude that we have to solve
\[\sum_{i=1}^j \pi_i\frac{1}{f_j-\kappa_i} = \sum_{i=j+1}^d \pi_i\frac{1}{\kappa_i-f_j}.\]
Observe that the left-hand side is $\infty$ for $f_j=\kappa_j$ and is decreasing in $[\kappa_j,\kappa_{j+1})$, whereas the right-hand side is increasing in the same interval and is $\infty$ for $f_j=\kappa_{j+1}$. 
This implies that there is a unique solution to the equation, which we simply call $f_j$. 
We have thus determined $\theta_j'(0) = f_j$, for $j=1,\ldots,d-1$.

\noindent
We proceed by computing $\theta_d'(0)$. 
From the relation $q^{-1} =\Psi_s(\theta(s))$ it follows by implicit differentiation to $s$, for any $j\in\{1,\ldots,d\}$,
with $\xi_{ij}(s):=q-\lambda_i(\ee^s-1)-\mu_i(\ee^{-s}-1)+\theta_j(s)$,
\[0=\sum_{i=1}^d\pi_i \frac{-\lambda_i \ee^s+\mu_i \ee^{-s}+\theta_j'(s)}{\xi_{ij}(s)^2}.\]
Inserting $s=0$ directly yields that $\theta_d'(0) = -{\boldsymbol\pi}{\boldsymbol\mu}\,(1-\rho)$, which in the sequel we will write as 
$-(1-\rho)\, {\mathbb E}{\mathcal M}$. 

\item[$\circ$] Differentiating once more, we arrive at
\begin{align*}
0=2\sum_{i=1}^d\pi_i \frac{\left(-\lambda_i \ee^s+\mu_i \ee^{-s}+\theta_j'(s)\right)^2}{\xi_{ij}(s)^3}
-\sum_{i=1}^d\pi_i \frac{-\lambda_i \ee^s-\mu_i \ee^{-s}+\theta_j''(s)}{\xi_{ij}(s)^2}.
\end{align*}
Again plugging in $s=0$, we obtain in self-evident notation (with e.g.\  $ {\mathbb E}\Lambda:={\boldsymbol\pi} {\boldsymbol \lambda}$) that
\[\theta''_d(0) = {\mathbb E}\Lambda+  {\mathbb E}{\mathcal M}+\frac{2\, {\mathbb E}(\Lambda-{\mathcal M)}^2}{q}.
\] 
The values of $\theta_1''(0),\ldots,\theta_{d-1}''(0)$ are not relevant in our analysis, as will turn out below. 
\end{itemize}
Using that a determinant is the product of the eigenvalues, we have that $e(s):={\rm det}\,E(s)=\theta_1(s)\cdots\theta_d(s).$
Noting  that $e(0)=0$, we wish to expand $e(s)$ as $e'(0) s+ \tfrac{1}{2}e''(0) s^2 +O(s^3)$ as $s\to 0$. 
Using that $\theta_d(0)=0$ and $\theta_j(0)=-q$ for $j=1,\dots,d-1$ we find, by applying the standard rules for differentiation of products,  that
\[e'(0) = \theta_d'(0) \prod_{j=1}^{d-1}\theta_j(0) = -(1-\rho)\, {\mathbb E}{\mathcal M}\cdot (-q)^{d-1}.\]
Likewise, for the second derivative we find
\begin{align*}e''(0) &= \theta_d''(0) \prod_{j\not= d}\theta_j(0) +\theta_d'(0)\sum_{j=1}^{d-1} \theta_j'(0)\prod_{k\not=j,d} \theta_k(0)\\
&=\theta''_d(0)\cdot(-q)^{d-1}-
(1-\rho)\, {\mathbb E}{\mathcal M}\cdot\sum_{j=1}^{d-1}f_j \cdot (-q)^{d-2}.
\end{align*}
The next step is to use the above findings to {expand the determinant $e(s)$ of $E(s)=A(\ee^{s})$ around $s=0$ and evaluated it in $-(1-\rho)s$}.
Similar to Eqn.\ \eqref{ch7:expand}, we obtain that, for any given $s$,
\begin{align*}e(-(1-\rho)s) &= -e'(0)\,(1-\rho)s + \tfrac{1}{2} e''(0)\,(1-\rho)^2s^2+O\left((1-\rho)^3\right)\\
&= (1-\rho)^2 \, {\mathbb E}{\mathcal M}\cdot (-q)^{d-1} s+\tfrac{1}{2} \,(1-\rho)^2\theta''_d(0)\cdot(-q)^{d-1} s^2+ O\left((1-\rho)^3\right).
\end{align*}
Next we can, in a fully analogous fashion, do the same for the determinants $e_i(s)$ of $E_i(s):=A_i(\ee^{s})$, and sum these over $i$. 
We thus obtain, by dividing the numerator and denominator by $(1-\rho)^2 s$, that
\[\lim_{\rho \uparrow 1} {\mathbb E}[\ee^{-s(1-\rho)Q}] = \frac{1}{1+\delta s},\]
with $\delta :=  {\theta''_d(0)}/({2\,  {\mathbb E}{\mathcal M}})$; in the limiting regime ($\rho \uparrow 1$), {we find that $\theta''_d(0)$ converges to}
\[{\lim_{\rho\uparrow 1} \theta''_d(0) }=\sigma^2_{\Lambda,{\mathcal M}}:=\big({\mathbb V}{\rm ar}\,\Lambda-2\,{\mathbb C}{\rm ov}(\Lambda,{\mathcal M})+{\mathbb V}{\rm ar}\,{\mathcal M}\big)\frac{2}{q}+2\,  {\mathbb E}{\mathcal M}.\] 

\begin{theorem} \label{ch7:TH2} 
Consider the resampled M/M/1 queue.
As $\rho\uparrow 1$, we have that $(1-\rho)Q$ converges to an exponentially distributed random variable with mean $\sigma^2_{\Lambda,{\mathcal M}}/(2 \,  {\mathbb E}{\mathcal M})$. 
\end{theorem}

{Note that in this special case we are able to fully determine the term $\alpha^\circ=\lim_{\rho \uparrow 1}\alpha''(1)$ that appears in Thm.\ \ref{ch7:TH1}; it turns out that $\lim_{\rho \uparrow 1}\alpha''(1)/\tau$ equals $\sigma^2_{\Lambda,{\mathcal M}}$ and $\tau = (-q)^{d-1}$.}

\vb

\noindent
{\it Numerical example: convergence to heavy-traffic limiting distribution}.

We conclude this subsection with a numerical example, assessing the accuracy of an approximation based on the above theorem. 
In our experiment we consider a two-dimensional background process in which the arrival rate is resampled at exponentially distributed epochs with mean $q^{-1}$: it is $\lambda$ with probability $\pi_1=\frac12$ and $0$ else. 
Observe that this means that the uninterrupted time the arrival rate is $\lambda$ (0, respectively) is exponentially distributed with parameter $q/2$. 
The service rate $\mu=1$ we hold fixed. 
In this setting the load is $\rho=\lambda/2$. 
We compute the distribution of the scaled queue length, and compare it with its heavy-traffic limit as $\rho\uparrow 1$ (or, alternatively, $\lambda \uparrow 2)$.
 
As can be found in e.g.\ \cite{Neuts}, with $\zeta_{ki}:= {\mathbb P}(Q=k,J=i)$, the stationary queue length has a matrix-geometric distribution:
\[
( \zeta_{k1} \, \zeta_{k2})
= ( \zeta_{01} \, \zeta_{02})\, R^k,
\]
where the matrix $R$ is the minimal nonnegative solution of the nonlinear matrix equation 
\[
\left(\begin{array}{cc} 
\lambda&0\\
0&0
\end{array}\right) + R
\left(\begin{array}{cc} 
-q/2-\lambda-1&q/2\\
q/2&-q/2-1
\end{array}\right) +R^2 \left(\begin{array}{cc} 
1&0\\
0&1 
\end{array}\right)={\boldsymbol 0},
\]
with ${\boldsymbol 0}$ denoting a $2\times 2$ all-zeros matrix. 
The probabilities $\zeta_{0i}$, $i=1,2$ follow from 
\[
(\zeta_{01} \, \zeta_{02}) 
\left( \left(\begin{array}{cc} -q/2-\lambda&q/2\\
q/2&-q/2
\end{array}\right) +R \left(\begin{array}{cc} 1&0\\
0&1 \end{array}\right) \right) 
= (0 \, 0),
\]
in combination with the normalization $(\zeta_{01} \, \zeta_{02}) (I_2-R)^{-1}{\boldsymbol 1}=1$, where ${\boldsymbol 1} = (1 \, 1)^\top$ and $I_2$ is the two-dimensional identity matrix.
The entries of $R$ and $(\zeta_{01} \, \zeta_{02})$ can be found numerically.

\begin{table}[ht!]
\centering
\begin{tabular}{|r|ccc|}
\hline
$\lambda = 1.8$ & $q = 0.5$ & $q = 1$ & $q = 2$ \\
\hline
$R$
& $\left(\begin{array}{cc} 0.963
& 0.837
\\ 0 
& 0 
\end{array}\right)$ 
& $\left(\begin{array}{cc} 0.946
& 0.854
\\ 0 
& 0 
\end{array}\right)$ 
& $\left(\begin{array}{cc} 0.930
&0.870
\\ 0 
& 0 
\end{array}\right)$ 
\\ 
$(\zeta_{01} \, \zeta_{02})$
& $\left(\begin{array}{cc} 0.0187
& 0.0813
\end{array}\right)$ 
& $\left(\begin{array}{cc} 0.0270 
& 0.0730
\end{array}\right)$ 
& $\left(\begin{array}{cc} 0.0349
& 0.0652
\end{array}\right)$ 
\\
\hline
$\lambda = 1.9$ & $q = 0.5$ & $q = 1$ & $q = 2$ \\ 
\hline
$R$
& $\left(\begin{array}{cc} 0.982
& 0.918
\\ 0 
& 0 
\end{array}\right)$ 
& $\left(\begin{array}{cc} 0.974
& 0.926
\\ 0 
& 0 
\end{array}\right)$ 
& $\left(\begin{array}{cc} 0.966
& 0.934
\\ 0 
& 0 
\end{array}\right)$ 
\\ 
$(\zeta_{01} \, \zeta_{02})$
& $\left(\begin{array}{cc} 0.0088
& 0.0412
\end{array}\right)$ 
& $\left(\begin{array}{cc} 0.0130
& 0.0370
\end{array}\right)$ 
& $\left(\begin{array}{cc} 0.0171
& 0.0330
\end{array}\right)$ 
\\
\hline
\end{tabular}%
\caption{$R$ and $(\zeta_{01} \, \zeta_{02})$ for different values of $\lambda$ and $q$.}
\label{table:numerics3}
\end{table}

In Table \ref{table:numerics3}, $R$ and $(\zeta_{01} \, \zeta_{02})$ are given for different values of $\lambda$ and $q$.
We consider two cases where the system approaches heavy traffic, namely $\lambda = 1.8$ and $\lambda = 1.9$. 
The corresponding loads are $\rho = 0.9$ and $\rho = 0.95$, respectively.
Note that $R_{i,j}$ describes the probability that the background state is $j$ when level $n+1$ is first hit given that we start in $(n,i)$.
Observe that in this light it makes sense that $R_{1,1}$ decreases when $q$ gets larger (faster resampling), while $R_{1,2}$ increases.
Also observe that the probability of an empty system $\zeta_{01} + \zeta_{02}$ equals $1-\rho$, as expected.

\begin{figure}[htb!]
	\begin{minipage}[t]{.3\textwidth}
		\centering
		\includegraphics[width=\textwidth]{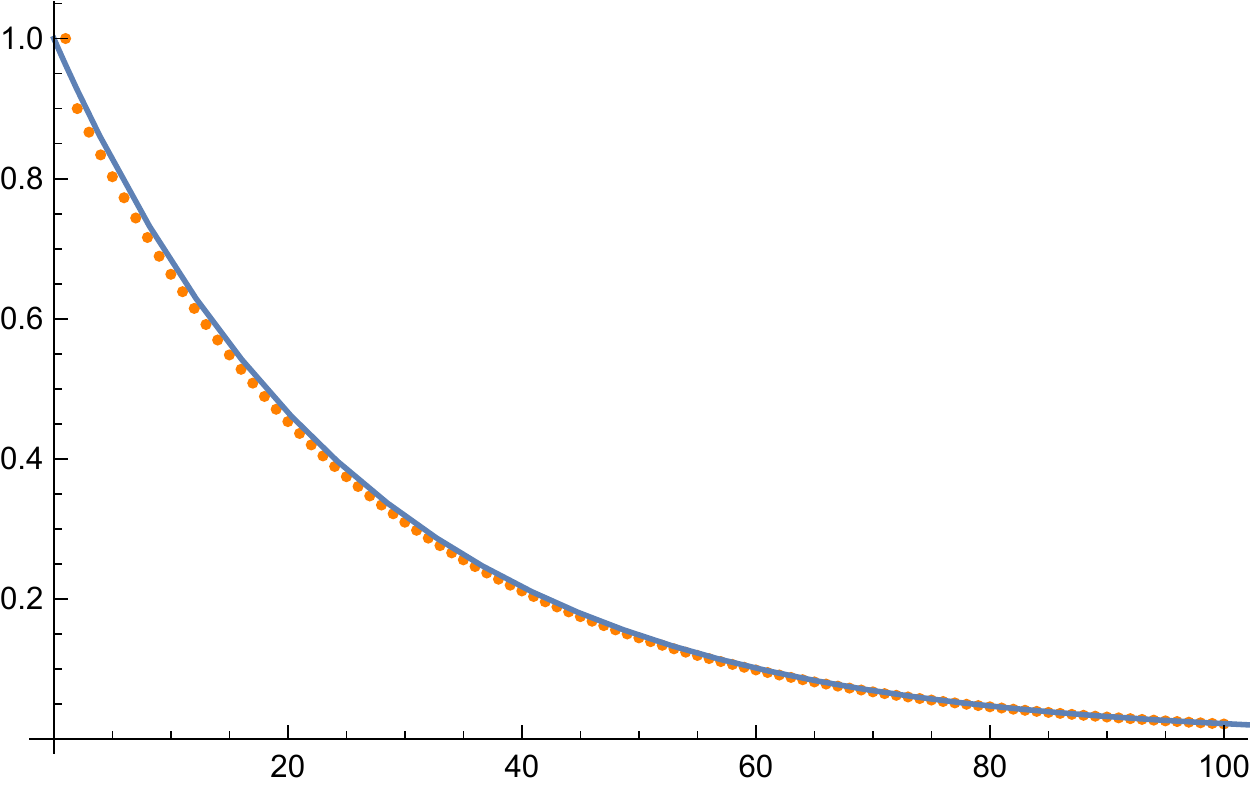} 
	\end{minipage}
	\hfill
	\begin{minipage}[t]{.3\textwidth}
		\centering
		\includegraphics[width=\textwidth]{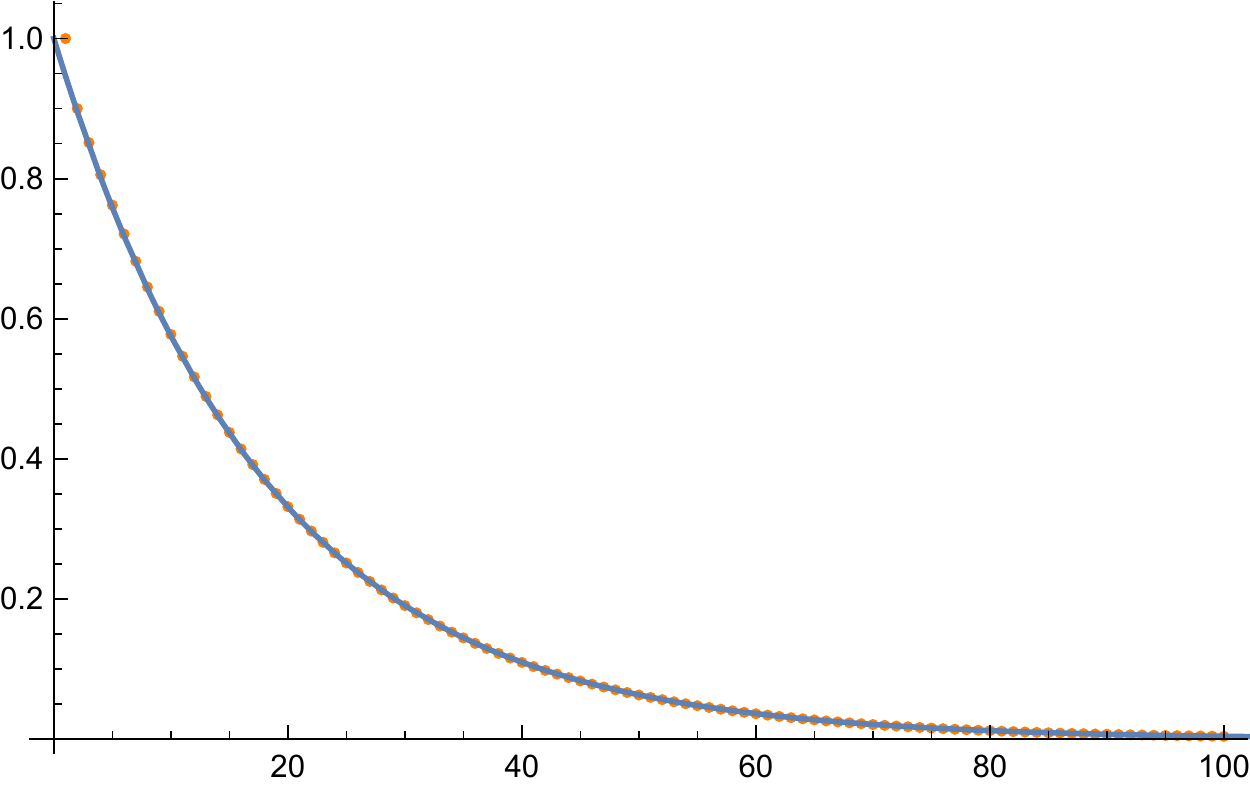} 
	\end{minipage}
	\hfill
	\begin{minipage}[t]{.3\textwidth}
		\centering
		\includegraphics[width=\textwidth]{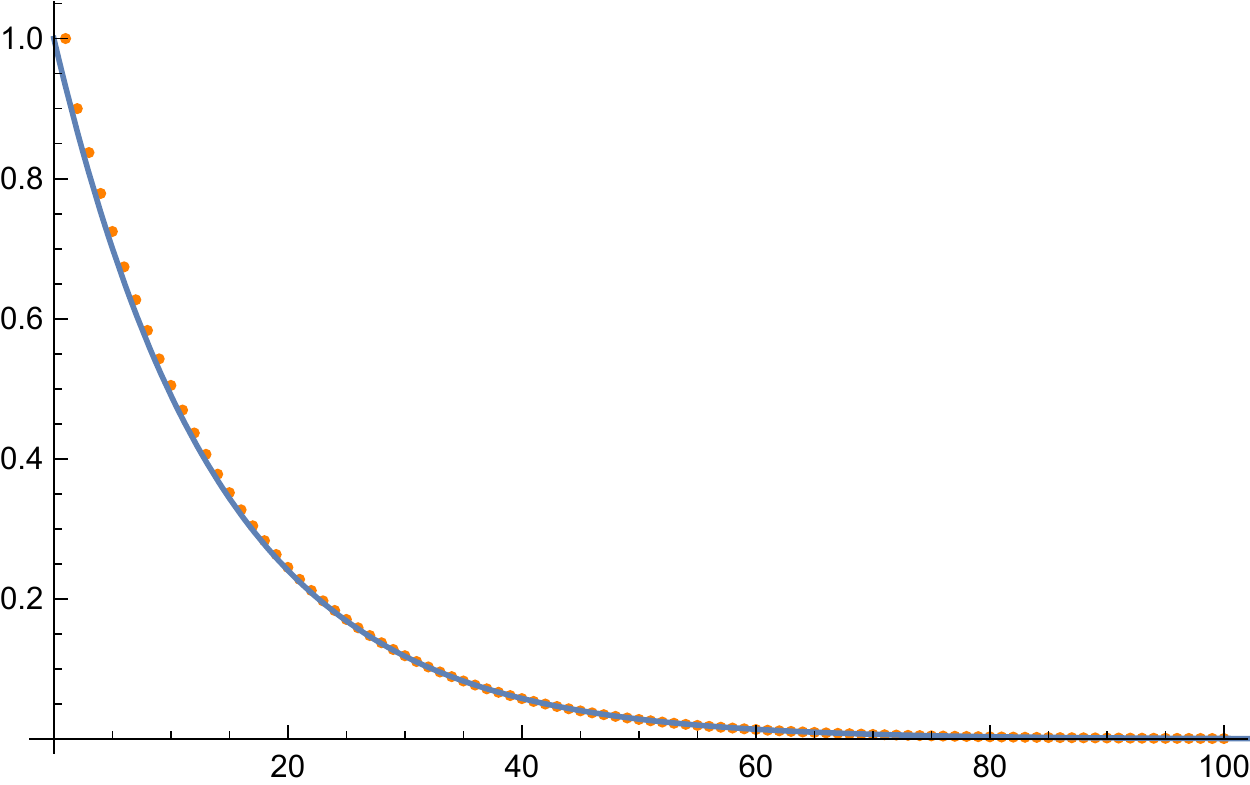} 
	\end{minipage}
	\hfill
	\begin{minipage}[t]{.3\textwidth}
		\centering
		\includegraphics[width=\textwidth]{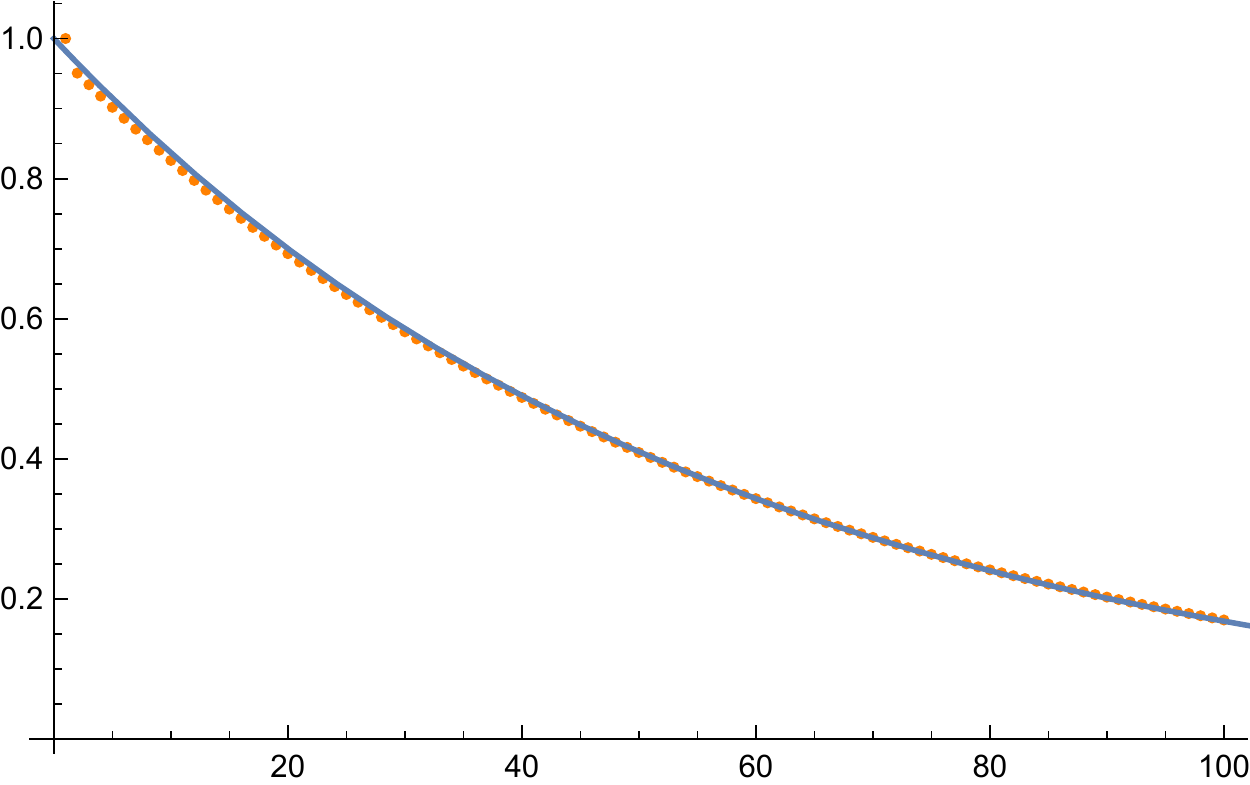} 
	\end{minipage}
	\hfill
	\begin{minipage}[t]{.3\textwidth}
		\centering
		\includegraphics[width=\textwidth]{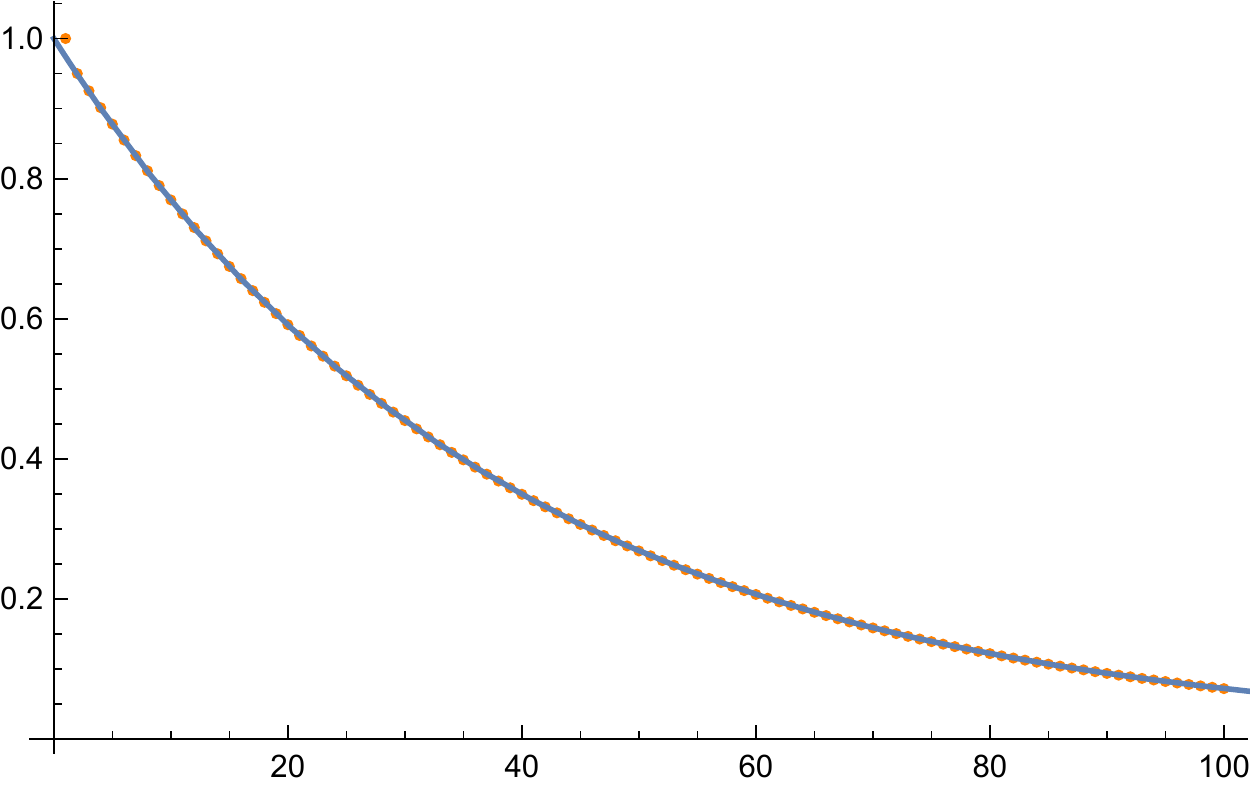}
	\end{minipage}
	\hfill
	\begin{minipage}[t]{.3\textwidth}
		\centering
		\includegraphics[width=\textwidth]{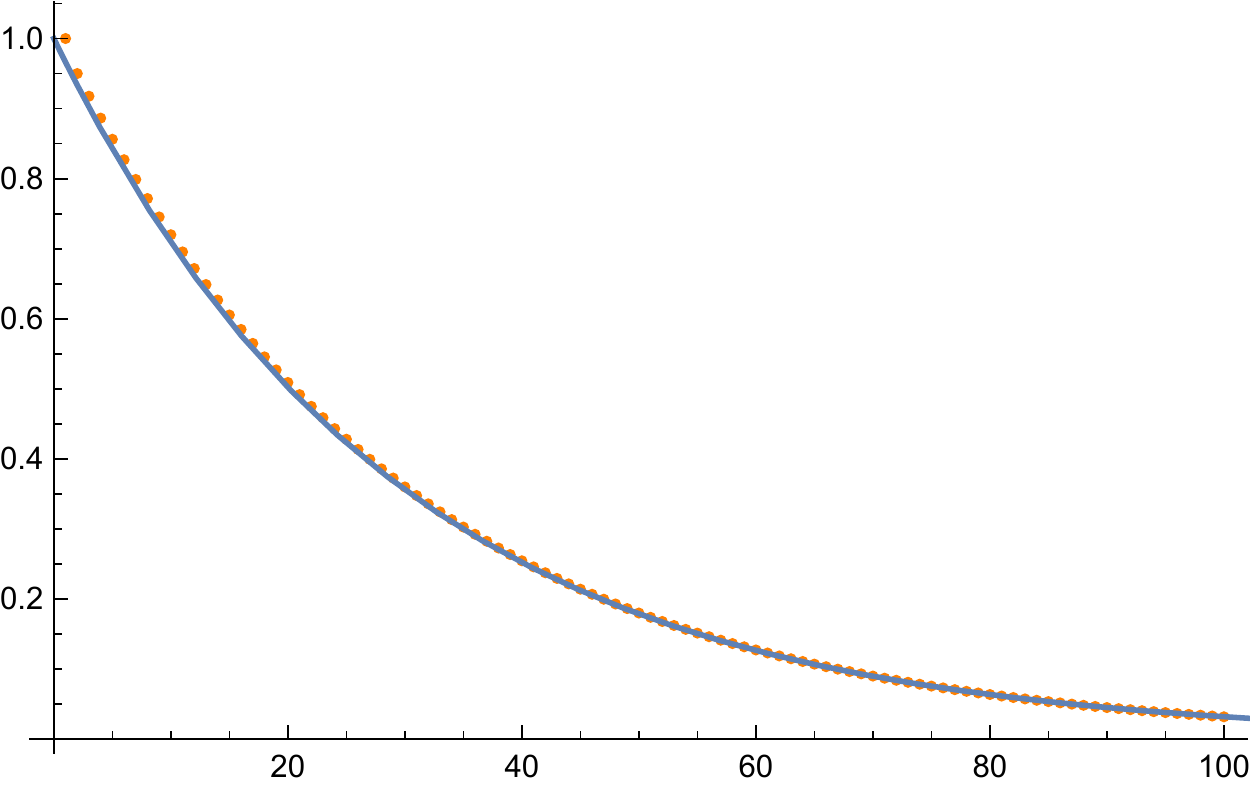} 
	\end{minipage}
	\hfill
\caption{ Plots of tail probabilities: exact solution (\textit{orange dotted line}) vs.\ exponential limiting distribution of the heavy-traffic result in Thm.\ \ref{ch7:TH2} (\textit{blue line}). 
The leftmost plots correspond to $q = 0.5$, the middle ones to $q=1$, and the rightmost ones to $q=2$; the top plots correspond to $\lambda = 1.8$, the bottom ones to $\lambda =1.9.$
\label{ch7:plotsHT}}
\end{figure}

Using $R$ and $(\zeta_{01} \, \zeta_{02})$ as given in Table \ref{table:numerics3}, we determine the exact tail probability ${\mathbb P}(Q > n)$ of the queue length, given by $1-\sum_{k=1}^{n-1} (\zeta_{k1} + \zeta_{k2})$; see the \textit{orange dotted line} in the plots in Figure \ref{ch7:plotsHT}. 
The numerical/exact results based on the above exact approach are to be compared with the limiting result of Theorem \ref{ch7:TH2}. To this end, observe that $ {\mathbb V}{\rm ar}\,\Lambda = \lambda^2/2-(\lambda/2)^2= \lambda^2/4,$ so that
\[\sigma^2_{\Lambda,{\mathcal M}}=\frac{\lambda^2}{2q}+2.\]
Hence, as $\lambda\uparrow 2$, $(1-\lambda/2)Q$ is approximately exponentially distributed with mean
$\lambda^2/(4q)+1$.
The \textit{blue line} in the plots in Figure \ref{ch7:plotsHT} depicts this exponential tail probability.

From the plots we conclude that, especially for the larger values of $q$, the approximation of the queue length distribution by the heavy-traffic limiting distribution is highly accurate for the considered values of $\rho$.
Note that the approximation gets near-perfect for any value of $q$ as we increase the load further.

\section{Exogenously triggered resampling, uncountable support}
\label{sec:exoperiodic}
The model considered in this section is an M/M/1 queue in which the arrival rate and service rate are resampled, but now the support of $(\Lambda,{\mathcal M})$ can be uncountable. 
Another difference with the previous section is that the resampling does not necessarily take place at Poisson epochs.
\subsection{Model description}
The inter-sample intervals are i.i.d.\ and distributed according to a non-negative random variable $\xi$ with density $f_\xi(\cdot)$.
We assume that at time $0$, the time until the next resampling (say $\bar\xi$) has the residual lifetime distribution, i.e., it has the density
\[f_{\bar\xi}(x):=\frac{{\mathbb P}(\xi>x)}{{\mathbb E}\xi}=\frac{1}{{\mathbb E}\xi}\int_x^\infty f_\xi(y)\,{\rm d}y;\] in other words, the resampling process is in stationarity. Above we tacitly assumed that ${{\mathbb E}\xi}<\infty$; in the sequel we also need that ${{\mathbb E}[\xi^2]}<\infty$.

We proceed by describing the cumulative input process $A(\cdot)$ and cumulative potential service process $S(\cdot)$. 
At any resample epoch the new arrival rate and service rate are sampled; these bivariate random quantities are i.i.d.\ and distributed as a componentwise non-negative two-dimensional random vector $(\Lambda,{\mathcal M})$.
In addition, we let $\Lambda(t)$ and ${\mathcal M}(t)$ be the arrival rate and service rate, respectively, that apply at time $t$. Throughout this section we assume that both ${\mathbb V}{\rm ar}\,\Lambda$ and ${\mathbb V}{\rm ar}\,{\mathcal M}$ are finite.

\subsection{Arrival and potential service process} \label{arpos}
As the resampling process starts in stationarity, it is clear that ${\mathbb E}A(t) = t\cdot{\mathbb E}\Lambda$ and ${\mathbb E}S(t) = t\cdot{\mathbb E}{\mathcal M}$. We therefore shift our attention to computing the variances and the covariance, relying on the law of total variance. 
The starting point is the representation of $A(t)$ and $S(t)$ in terms of Poisson processes with random parameter:
\[A(t) = {\rm Pois}\left(\int_0^t \Lambda(s)\,{\rm d}s\right),\:\:\:\:S(t) = {\rm Pois}\left(\int_0^t {\mathcal M}(s)\,{\rm d}s\right).\]
We show how the calculation for ${\mathbb V}{\rm ar}\,A(t)$ is done; the calculation for ${\mathbb V}{\rm ar}\,S(t)$ can be done fully analogously. The law of total variance entails, with $\Lambda(\cdot)\equiv (\Lambda(s))_{s\in[0,t]}$,
\begin{equation}
\label{dec}v_A(t):={\mathbb V}{\rm ar}\,A(t) =   {\mathbb E}[{\mathbb V}{\rm ar}(A(t)\,|\,\Lambda(\cdot))] +{\mathbb V}{\rm ar}\left({\mathbb E}[A(t)\,|\,\Lambda(\cdot)]\right) .\end{equation}
Clearly,
\[{\mathbb E}[A(t)\,|\,\Lambda(\cdot)] = {\mathbb V}{\rm ar}(A(t)\,|\,\Lambda(\cdot)) = \int_0^t \Lambda(s)\,{\rm d}s.\]
The first term in (\ref{dec}) therefore equals
\[  {\mathbb E}[{\mathbb V}{\rm ar}(A(t)\,|\,\Lambda(\cdot))] ={\mathbb E}\left[\int_0^t \Lambda(s)\,{\rm d}s\right] = \int_0^t {\mathbb E}\Lambda(s)\,{\rm d}s =  {\mathbb E}A(t) = t\cdot{\mathbb E}\Lambda.\]
Now focus on the second term in (\ref{dec}). Using standard properties, we obtain
\[{\mathbb V}{\rm ar}\left({\mathbb E}[A(t)\,|\,\Lambda(\cdot)]\right) = 2\int_0^t\int_0^s {\mathbb C}{\rm ov}\left(\Lambda(r),\Lambda(s)\right)\,{\rm d}r\,{\rm d}s.\] 
The crucial insight is that ${\mathbb C}{\rm ov}(\Lambda(r),\Lambda(s))$ equals ${\mathbb V}{\rm ar}\,\Lambda$ if $r$ and $s$ are in the same resampling interval, and $0$ else. 
With $\bar\xi$ denoting the stationary residual lifetime pertaining to $\xi$, we thus obtain
\[{\mathbb V}{\rm ar}\left({\mathbb E}[A(t)\,|\,\Lambda(\cdot)]\right) = 2\,{\mathbb V}{\rm ar}\,\Lambda\int_0^t\int_0^s {\mathbb P}(\bar\xi> s-r)\,{\rm d}r\,{\rm d}s.\] 
Notice that, by a standard calculation, with $g(u):={\mathbb P}(\xi> u)$,
\begin{align*}
{\mathbb E}\xi&\int_0^t\int_0^s {\mathbb P}(\bar\xi> s-r)\,{\rm d}r\,{\rm d}s=\int_0^t\int_0^s \int_{s-r}^\infty g(u)\,{\rm d}u\,{\rm d}r\,{\rm d}s
=\int_0^t\int_0^{t-r} \int_{y}^{\infty}  g(u)\,{\rm d}u\,{\rm d}y\,{\rm d}r\\
&=\int_0^{\infty}\int_0^{t} \int_{0}^{\min\{u,t-r\}} \hspace*{-8pt} g(u)\,{\rm d}y\,{\rm d}r\,{\rm d}u\\
&=\left(\int_0^t\left(\int_0^{t-u} \int_{0}^{u}  g(u)\,{\rm d}y\,{\rm d}r + \int_{t-u}^t \int_{0}^{t-r}  g(u)\,{\rm d}y\,{\rm d}r\right)\,{\rm d}u + \int_t^{\infty}\int_0^{t} \int_{0}^{t-r}  g(u)\,{\rm d}y\,{\rm d}r\,{\rm d}u \right)\\
&=\left(\int_0^t(tu-\tfrac{1}{2}u^2) \, g(u)\ {\rm d}u+\tfrac{1}{2}t^2 \int_t^\infty   g(u)\, {\rm d}u\right).
\end{align*}
We have thus found the following expression (and its counterpart for the service process); here $v_S(t):= {\mathbb V}{\rm ar}\,S(t) $.
\begin{proposition} The variances of the cumulative arrival and service process are given by, for $t\geqslant 0$,
\[v_A(t) = t\cdot{\mathbb E}\Lambda+2\frac{{\mathbb V}{\rm ar}\,\Lambda}{{\mathbb E}\xi}\left(\int_0^t(tu-\tfrac{1}{2}u^2) \,{\mathbb P}(\xi> u)\ {\rm d}u+\tfrac{1}{2}t^2 \int_t^\infty  {\mathbb P}(\xi> u)\, {\rm d}u\right),\]
and
\[v_S(t) = t\cdot{\mathbb E}{\mathcal M}+2\frac{{\mathbb V}{\rm ar}({\mathcal M})}{{\mathbb E}\xi}\left(\int_0^t(tu-\tfrac{1}{2}u^2) \,{\mathbb P}(\xi> u)\ {\rm d}u+\tfrac{1}{2}t^2 \int_t^\infty  {\mathbb P}(\xi> u)\, {\rm d}u\right).
\]
\end{proposition}
\begin{example}
It can be checked that for the special case that $\xi$ is exponentially distributed with parameter $\delta>0$ (implying that $\bar \xi$ is exponentially distributed with parameter $\delta$ as well),
\begin{equation}
\int_0^t\int_0^s {\mathbb P}(\bar\xi> s-r)\,{\rm d}r\,{\rm d}s =\frac{t}{\delta} -\frac{1}{\delta^2}(1-\ee^{-\delta t}).
\end{equation}
As a consequence
\[v_A(t) =  t\cdot{\mathbb E}\Lambda+2\,{{\mathbb V}{\rm ar}\,\Lambda}\left(
\frac{t}{\delta} -\frac{1}{\delta^2}(1-\ee^{-\delta t})
\right),\]
for $t\geqslant 0.$ \hfill$\Box$
\end{example}
The next step is to evaluate ${\mathbb C}{\rm ov}(A(t),S(t))$ by the law of total covariance. We have, with $P(\cdot)\equiv(\Lambda(s),{\mathcal M}(s))_{s\in[0,t]}$ denoting the parameter process in the interval $s\in[0,t]$, 
\[c_{A,S}(t):={\mathbb C}{\rm ov}(A(t),S(t)) = {\mathbb E}[{\mathbb C}{\rm ov}(A(t),S(t)\,|\,P(\cdot))]+ {\mathbb C}{\rm ov}({\mathbb E}[A(t)\,|\,P(\cdot)] , {\mathbb E}[S(t)\,|\,P(\cdot)]).\]
The first term on the right-hand side is clearly 0: conditioned on $P(\cdot)$, the arrival and service processes are independent. The second term on the right-hand side can be dealt with as before. We thus obtain the following result. \begin{proposition} The covariance between the cumulative arrival and service process is given by, for $t\geqslant 0$,
\[c_{A,S}(t) = 
2\frac{{\mathbb C}{\rm ov}(\Lambda,{\mathcal M})}{{\mathbb E}\xi}\left(\int_0^t(tu-\tfrac{1}{2}u^2) \,{\mathbb P}(\xi> u)\ {\rm d}u+\tfrac{1}{2}t^2 \int_t^\infty  {\mathbb P}(\xi> u)\, {\rm d}u\right).
\]
\end{proposition}
Our next objective is to show that $v_A(t)$, $v_S(t)$, and $c_{A,S}(t)$ behave essentially linear as $t\to\infty$; this finding will later play a role in the derivation of functional limit theorems. 
We start by considering $v_A(t)/t$ as $t\to\infty.$ By L'H\^opital's theorem, this limit equals
\[\lim_{t\to\infty} v_A'(t)={\mathbb E}\Lambda+\lim_{t\to\infty} 2\frac{{\mathbb V}{\rm ar}\,\Lambda}{{\mathbb E}\xi} \int_0^\infty \min\{u,t\} {\mathbb P}(\xi>u)\,{\rm d}u.\]
By a straightforward calculation (applying dominated convergence and integration by parts), we obtain 
\[v_A:=\lim_{t\to\infty} \frac{v_A(t)}{t} = {\mathbb E}\Lambda+{{\mathbb V}{\rm ar}\,\Lambda}\frac{{\mathbb E}[\xi^2]}{{\mathbb E}\xi}.\] The limit of $v_S(t)/t$ as $t\to\infty$ (which we call $v_S$) and the limit of $c_{A,S}(t)/t$ as $t\to\infty$ (which we call $c_{A,S}$) can be determined in the same fashion.
We state this result as a corollary.
\begin{corollary}\label{C1}The constants $v_A$, $v_S$, and $c_{A,S}$ are given by
\[ v_A= {\mathbb E}\Lambda+{{\mathbb V}{\rm ar}\,\Lambda}\,\frac{{\mathbb E}[\xi^2]}{{\mathbb E}\xi},\:\:\:\:v_S= {\mathbb E}\mathcal M+{{\mathbb V}{\rm ar}\,\mathcal M}\,\frac{{\mathbb E}[\xi^2]}{{\mathbb E}\xi},\:\:\:\:c_{A,S}={{\mathbb C}{\rm ov}(\Lambda,\mathcal M)}\,\frac{{\mathbb E}[\xi^2]}{{\mathbb E}\xi}.\]
\end{corollary}
An alternative derivation of these expressions can be found in Appendix \ref{APP}; the methodology presented there is particularly useful, as it also facilitates the derivation of higher moments in a relatively straightforward manner. 

\subsection{Weak convergence to reflected Brownian motion} \label{sec:WCBM}
After having studied some properties of the arrival and potential service process, we now shift our attention to the queueing process. Our objective is to establish in the heavy-traffic regime, under a time-scaling, weak convergence of the queueing process to reflected Brownian motion.

Assuming the system starts empty at time $0$, the number of customers in the queue at time $t$ can be written as
\[\sup_{0\leqslant s\leqslant t} \big( A(t)- A(s) -(S(t) - S(s))\big);\]
this representation involving the potential service process applies due to the fact that the service times are exponential (albeit with some random value), cf.\ the remark in \cite[bottom of p.\ 290]{WHI}. 

In this subsection
we impose the parameterization ${\mathbb E}\Lambda = \rho \ {\mathbb E}{\mathcal M}$, and consider the regime $\rho\uparrow 1$ (cf.\ Remark \ref{rem:scaling} for more details). 
Clearly, the queueing process, which we denote by $(Q_\rho(t))_{t\geqslant 0}$ to stress the dependence on $\rho$, blows up as $\rho\uparrow 1.$ 
However, after appropriately rescaling time and space, one obtains a non-trivial limiting process. 
More specifically, we will study the behavior of
\begin{equation}\label{Qrho}\left((1-\rho) Q_\rho\left(\frac{t}{(1-\rho)^2}\right)\right)_{t\geqslant 0}\end{equation}
for $\rho\uparrow 1$; observe that time is stretched by a factor $(1-\rho)^2$ where space is compressed by a factor $1-\rho$. We do so by showing that the process $(B_\rho(t))_{t\geqslant 0}$, with \[B_\rho(t) := (1-\rho) \,A_\rho\left(\frac{t}{(1-\rho)^2}\right) -(1-\rho) \,S\left(\frac{t}{(1-\rho)^2}\right), \]
converges weakly  to a Brownian motion; we write $A_\rho(\cdot)$, with subscript $\rho$, to make visible that we choose ${\mathbb E}\Lambda$ equal to $\rho\ {\mathbb E}{\mathcal M}.$
By `continuous mapping' this convergence then also yields that (\ref{Qrho}) converges weakly to reflected Brownian motion.

Let $P_1(\cdot)$ and $P_2(\cdot)$ be defined as two independent unit-rate Poisson processes. Then, with $f_\rho:=1/(1-\rho)^2$,
\begin{equation}
\label{difp}B_\rho(t) = (1-\rho) P_1\left(\int_0^{tf_\rho} \Lambda_\rho(s)\,{\rm d}s\right) - (1-\rho) P_2\left(\int_0^{tf_\rho} {\mathcal M}(s)\,{\rm d}s\right),\end{equation}
where the subscript $\rho$ has been added to $\Lambda(\cdot)$ to indicate the dependence on $\rho$. 

In the first place,
\[\lim_{\rho\uparrow 1}{\mathbb E}B_\rho(t) = \lim_{\rho\uparrow 1}\frac{t}{1-\rho} \big({\mathbb E}\Lambda-{\mathbb E}{\mathcal M}\big) = -t\, {\mathbb E}{\mathcal M}. \]
Also, representation (\ref{difp}) implies that, for some martingale $K_\rho(\cdot)$,
\[ {\rm d}B_\rho(t) = (1-\rho)f_\rho\, \Lambda_\rho(tf_\rho)\,{\rm d}t - (1-\rho)f_\rho\, {\mathcal M}(tf_\rho)\,{\rm d}t + {\rm d}K_\rho(t).\]
Using similar computations as before,
\begin{align*}
\lefteqn{\lim_{\rho\uparrow 1}{\mathbb V}{\rm ar}\left((1-\rho) \int_0^{tf_\rho} \Lambda_\rho(s)\,{\rm d}s - (1-\rho) \int_0^{tf_\rho} {\mathcal M}(s)\,{\rm d}s\right)}\hspace{2cm}\\
=&\: \lim_{\rho\uparrow 1}(1-\rho)^2\,{\mathbb V}{\rm ar}\left( \int_0^{tf_\rho} \Lambda_\rho(s)\,{\rm d}s -  \int_0^{tf_\rho} {\mathcal M}(s)\,{\rm d}s\right)\\
=&\: \lim_{\rho\uparrow 1}(1-\rho)^2\,{\mathbb V}{\rm ar}(\Lambda-{\mathcal M})
\int_0^{tf_\rho} \int_0^s {\mathbb P}(\bar\xi>s-r){\rm d}r\,{\rm d}s
\\
=&\:t\,\left({\mathbb V}{\rm ar}\Lambda-2\, {\mathbb C}{\rm ov}(\Lambda,{\mathcal M})+{\mathbb V}{\rm ar}{\mathcal M} \right)\frac{{\mathbb E}[\xi^2]}{{\mathbb E}[\xi]}.
\end{align*}
In addition, with $\langle X(\cdot) \rangle_t$ denoting the quadratic variation process of $X(\cdot)$ at time $t$, using standard properties of pure jump processes,
\[\frac{\rm d}{{\rm d}t}\langle K_\rho(\cdot)\rangle_t =  (1-\rho)^2f_\rho\, \Lambda_\rho(tf_\rho)+ (1-\rho)^2f_\rho\, {\mathcal M}(tf_\rho)=\Lambda_\rho(tf_\rho)+  {\mathcal M}(tf_\rho),\]
so that, as $\rho\uparrow 1$,
\[ \langle K_\rho(\cdot)\rangle_t = \frac{1}{f_\rho}\left(\int_0^{tf_\rho} \Lambda_\rho(s)\,{\rm d}s+  \int_0^{tf_\rho} {\mathcal M}(s)\,{\rm d}s\right)\to 2t\cdot {\mathbb E}{\mathcal M} .\]
Combining the above, we conclude that, as $\rho\uparrow 1$, $B_\rho(\cdot)$ converges weakly to a process $B(\cdot)$ given by
\[B(t) = -t\, {\mathbb E}{\mathcal M} +\sigma_{\Lambda, {\mathcal M}}\cdot W(t),\]
with $W(\cdot)$ a standard Brownian motion and
\begin{equation} \label{varBM}
\sigma_{\Lambda, {\mathcal M}}^2:=\left({\mathbb V}{\rm ar}\Lambda-2\, {\mathbb C}{\rm ov}(\Lambda,{\mathcal M})+{\mathbb V}{\rm ar}{\mathcal M} \right)\frac{{\mathbb E}[\xi^2]}{{\mathbb E}[\xi]}+2\ {\mathbb E}{\mathcal M}.
\end{equation}
After applying the continuous mapping theorem, we have thus found the following result.
\begin{theorem} \label{RBM} As $\rho\uparrow 1$, with $Q(t):=\sup_{0\leqslant s\leqslant t}(B(t)-B(s))$,
\[\left((1-\rho) Q_\rho\left(\frac{t}{(1-\rho)^2}\right)\right)_{t\geqslant 0} \stackrel{\rm d}{\to} (Q(t))_{t\geqslant 0}.\]
\end{theorem}

\begin{remark}{\em 
In this remark we analyze how large $\sigma_{\Lambda, {\mathcal M}}^2$ can be. Given that both $\Lambda$ and ${\mathcal M}$ are non-negative, one is inclined to believe that one cannot achieve a correlation coefficient between $\Lambda$ and ${\mathcal M}$ with value $-1$. This is, however, not true, as follows from the following argument.

Consider two  non-negative random variables $X$ and $Y$, having (without losing any generality) unit mean. Let $X$ have a  given distribution, and define
\[Y:= \frac{1}{\psi(s)}\left(-\frac{1}{s}X
+1\right)1_{\{0\leqslant X\leqslant s\}},\]
where, with $p(s):={\mathbb P}(X\leqslant s)$,
\[\psi(s):={\mathbb E}\left[\left(-\frac{1}{s}X
+1\right)1_{\{0\leqslant X\leqslant s\}}\right]=-\frac{1}{s}{\mathbb E}\left[X\,
1_{\{0\leqslant X\leqslant s\}}\right]+p(s).
\]
It is immediately seen that $Y$ is indeed non-negative with mean equal to 1.

Now, under obvious mild regularity conditions that ensure the existence of the expectations involved, as $s\to\infty$,
\begin{align*}{\mathbb C}{\rm ov}(X,Y)=&\:{\mathbb E}[XY]-1 =\frac{1}{\psi(s)}{\mathbb E}\left[\left(-\frac{1}{s}X^2
+X\right)1_{\{0\leqslant X\leqslant s\}}\right] -1\\
=&\:\frac{-s^{-1} {\mathbb E}
\left[X^2\,1_{\{0\leqslant X\leqslant s\}}\right]+
(1+s^{-1}){\mathbb E}\left[X\,1_{\{0\leqslant X\leqslant s\}}\right]-
p(s)}
{-s^{-1}{\mathbb E}\left[X\,
1_{\{0\leqslant X\leqslant s\}}\right]+p(s)}\\
\sim&\:\frac{-{\mathbb V}{\rm ar}(X)-({\mathbb E}[X])^2+{\mathbb E}[X]}{s}=-\frac{{\mathbb V}{\rm ar}(X)}{s}.
\end{align*}
Likewise, again as $s\to\infty$,
\begin{align*}{\mathbb V}{\rm ar}(Y)=&\:{\mathbb E}[Y^2]-1 =\frac{1}{\psi^2(s)}{\mathbb E}\left[\left(\frac{1}{s^2}X^2
-\frac{2}{s}X+1\right)1_{\{0\leqslant X\leqslant s\}}\right] -1\\
=&\:\frac{s^{-2} {\mathbb E}
\left[X^2\,1_{\{0\leqslant X\leqslant s\}}\right]
-2s^{-1}{\mathbb E}\left[X\,1_{\{0\leqslant X\leqslant s\}}\right]+p(s)-
(-{s}^{-1}{\mathbb E}\left[X\,
1_{\{0\leqslant X\leqslant s\}}\right]+p(s))^2
}
{(-s^{-1}{\mathbb E}\left[X\,
1_{\{0\leqslant X\leqslant s\}}\right]+p(s))^2}\\
\sim&\:\frac{{\mathbb V}{\rm ar}(X)}{s^2}.
\end{align*}
Conclude that this choice of $Y$ (for a given $X$) yields a correlation coefficient $-1$ as $s$ grows large. 
This means that one can achieve a correlation coefficient arbitrarily close to $-1$. 

Applying this observation to the expression for the variance in $(\ref{varBM})$ for any $\varepsilon>0$ and any non-negative $\Lambda$, one can construct a non-negative ${\mathcal M}$ that is negatively correlated with $\Lambda$ in such a way that \[\sigma_{\Lambda, {\mathcal M}}^2 \geqslant \big(\sqrt{{\mathbb V}{\rm ar}\Lambda}+\sqrt{{\mathbb V}{\rm ar}{\mathcal M}} \big)^2\,\frac{{\mathbb E}[\xi^2]}{{\mathbb E}[\xi]}+2\ {\mathbb E}{\mathcal M}-\varepsilon.\]}
\end{remark}

\begin{remark}
{\em
Based on Theorem \ref{RBM} one would expect that the stationary number of customers under the heavy-traffic scaling would converge to the stationary version of $B(\cdot)$ reflected at 0, which has an exponential distribution with mean $\sigma^2_{\Lambda,{\mathcal M}}/(2\,\E{\mathcal M})$. Establishing such a result, however, would require interchanging two limits, namely $t\to\infty$ and $\rho\uparrow 1$. 
For specific models, an argumentation developed in \cite{SW} can be followed, but it is not clear how this could be applied in our setting.
}
\end{remark}

\section{Endogenously triggered resampling}
\label{sec:endo}
In this section we consider an M/G/1 queue with the special feature that at every service completion the arrival rate is resampled. We refer to this mechanism as `endogenously triggered resampling', as the resampling is not due to an exogenous, independently evolving process. 
\subsection{Model description}
In the model we consider, the service times are i.i.d.\ samples from some general non-negative distribution, distributed as a generic random variable $S$ with Laplace-Stieltjes transform $\sigma(\cdot)$. 
The distinguishing feature of the model we study in this section is that at every service completion the arrival rate is resampled from a general distribution with non-negative support; the sequence of arrival rates is assumed i.i.d., each of them being distributed  as the generic random variable $\Lambda$. 

To make sure the queueing system under study is stable, we have to assume ${\mathbb E}\Lambda \leqslant {\mathbb E}S$.
Later, when considering the heavy-traffic regime, we will assume the finiteness of the corresponding second moments. 

\subsection{Transform of the stationary number of customers}
Let $N_n$ be the number of customers arriving during the $n$-th service time, and $Q_n$ the stationary number of customers present at the $n$-th service completion. It is evident that $Q_{n+1} =_{\rm d} (Q_n-1)^+ + N_{n+1},$ with the two terms in the right-hand side of this distributional equality being independent. We obtain the following relation:
\begin{align*}
 {\mathbb E} z^{(Q_n-1)^+} &= \sum_{j=0}^\infty z^j \,{\mathbb P}((Q_n-1)^+=j) 
={\mathbb P}(Q_n=0) + \sum_{j=1}^\infty z^{j-1} {\mathbb P}(Q_n= j)\\
&={\mathbb P}(Q_n=0) + \frac{1}{z}\left(\sum_{j=0}^\infty z^j {\mathbb P}(Q_n= j) - {\mathbb P}(Q_n=0) \right)\\
&=\left(1-\frac{1}{z}\right){\mathbb P}(Q_n=0)  + \frac{ {\mathbb E}z^{Q_n}}{z}.
\end{align*}
With, for $n\in{\mathbb N}_0$, $\kappa_n(z):=  {\mathbb E}z^{Q_n}$ and $\nu(z):= {\mathbb E} z^{N_n}$ (which evidently does not depend on $n$), we thus obtain the recursion
\[\nu(z) \left(\left(1-\frac{1}{z}\right) {\mathbb P}(Q_n=0) +\frac{\kappa_n(z)}{z}\right) = \kappa_{n+1}(z).\]
Our next step is that we consider the transform of $\kappa_n(z)$ at a geometrically distributed time epoch, which, as it turns out, can be expressed in closed form. To this end, we multiply the above recursion by $(1-r)^nr$ (for some $r\in[0,1)$) and sum over $n\in{\mathbb N}_0$. We thus obtain
\[
\sum_{n=0}^\infty (1-r)^nr \,\nu(z) \left(\left(1-\frac{1}{z}\right) {\mathbb P}(Q_n=0) +\frac{\kappa_n(z)}{z}\right) = \sum_{n=0}^\infty (1-r)^nr\,\kappa_{n+1}(z);
\]
(technically, the geometric distribution is a {\it shifted} geometric distribution, having probability mass  at $0$). 
Our objective is to identify the double transform
\[K(r,z) := \sum_{n=0}^\infty (1-r)^nr \,\kappa_n(z),\]
assuming that we know the distribution of $Q_0$ (i.e., $\kappa_0(z)$ is known). Observing that ${\mathbb P}(Q_n=0)= \kappa_n(0)$, we thus obtain the identity
\[\nu(z)\left(1-\frac{1}{z}\right) K(r,0) + \frac{\nu(z)}{z}K(r,z) = \frac{1}{1-r}\big(K(r,z) - r\kappa_0(z)\big).\]
Isolating $K(r,z)$, we can express this double transform in terms of $K(r,0)$ and $\kappa_0(z)$: some elementary algebra yields
\[
K(r,z) = \frac{rz\,\kappa_0(z) - (1-r)(1-z)\nu(z)\,K(r,0)}{z-(1-r)\nu(z)}.
\]
The unknown function $K(\cdot\,,0)$ can be eliminated by observing that, for a given $r\in[0,1]$, each root of the denominator must correspond to a root of the numerator. Observe that, for a given $r\in(0,1)$, the function $z\mapsto (1-r) \nu(z)$ is increasing and convex, and attains values in $(0,1)$, which immediately implies that $z=(1-r)\nu(z)$ has a unique root $z_0(r)\in(0,1).$ The fact that this is a root of the numerator as well yields that
\[K(r,0) = \frac{rz_0(r)\,\kappa_0(z_0(r))}{(1-r)(1-z_0(r))\nu(z_0(r))}=\frac{r\,\kappa_0(z_0(r))}{1-z_0(r)} .\]
We thus arrive at the following result.
\begin{theorem} \label{thtr} For $r,z\in(0,1)$,
\[K(r,z) = \frac{rz\,\kappa_0(z)}{z-(1-r)\nu(z)} - \frac{(1-r)(1-z)\nu(z)}{z-(1-r)\nu(z)}\cdot \frac{r\,\kappa_0(z_0(r)) }{1-z_0(r)}.\]
\end{theorem}

In the same way, the stationary behavior can be dealt with; the object of study is the stationary number of customers $Q$. 
We have to impose ${\mathbb E}N<1$ to ensure the existence of a stationary distribution. 
It is evident that, under the stability assumption imposed,  $Q =_{\rm d} (Q-1)^+ + N,$ with the two terms in the right-hand side being independent. Then a standard argumentation, similar to the one used above, yields
\begin{align*}
 {\mathbb E} z^{(Q-1)^+} &= \sum_{j=0}^\infty z^j \,{\mathbb P}((Q-1)^+=j) 
={\mathbb P}(Q= 0) + \sum_{j=1}^\infty z^{j-1} {\mathbb P}(Q= j)\\
&={\mathbb P}(Q= 0)+ \frac{1}{z}\left(\sum_{j=0}^\infty z^j {\mathbb P}(Q= j) -  {\mathbb P}(Q= 0)\right)
=\left(1-\frac{1}{z}\right){\mathbb P}(Q= 0) + \frac{ {\mathbb E}z^Q}{z}.
\end{align*}
Upon combining the above,
with $\kappa(z) :=  {\mathbb E}z^Q$,
\[\kappa(z) = \nu(z)\left(\left(1-\frac{1}{z}\right){\mathbb P}(Q= 0) +\frac{\kappa(z)}{z}\right).\]
Solving $\kappa(z)$ from this equation directly yields (where one should recall that $\nu(z)\geqslant z$ for all $z\in[0,1]$, which is due to ${\mathbb P}(N=0)>0$ and $\nu'(1)={\mathbb E}N<1$, in combination with the fact that $\nu(\cdot)$ is convex on $[0,1]$)
\[\kappa(z) = \nu(z)\frac{1-z}{\nu(z)-z}\,{\mathbb P}(Q= 0).\]
Because of $\kappa(1) = 1$, we find by L'H\^opital's rule that ${\mathbb P}(Q= 0) = 1-\nu'(1) = 1 -{\mathbb E}N.$

The above reasoning has been frequently used to obtain the stationary distribution in the M/G/1 queue, but it is important to observe that it does not use that the interarrival times are i.i.d.\ random variables from the same exponential distribution. 
Indeed, it allows for the arrival rate to be resampled at every service completion. 
In the resampling model we study in this section, with $g(\cdot)$ denoting the density of $\Lambda$,
\begin{equation}
\label{nuz}
\nu(z) = \int_0^\infty \sigma\big(\lambda(1-z)\big) \,g(\lambda)\,{\rm d}\lambda= {\mathbb E}\sigma\big(\Lambda(1-z)\big).\end{equation}
We end up with the following result, using that ${\mathbb E}N = {\mathbb E}\Lambda\ {\mathbb E}S$.
\begin{theorem} Under the stability constraint ${\mathbb E}\Lambda\ {\mathbb E}S<1$, for $z\in(0,1)$,\begin{equation}\label{LSTkappa}\kappa(z) = \nu(z)\frac{1-z}{\nu(z)-z} \left(1 -{{\mathbb E}\Lambda}\,{{\mathbb E}S}\right)
,\end{equation}
with $\nu(\cdot)$ given by $(\ref{nuz})$.
\end{theorem}

\subsection{Heavy-traffic scaling limit} \label{sec:endoHT}
In this subsection we consider the regime in which $\rho:= {{\mathbb E}\Lambda}\,{{\mathbb E}S}$ goes to $1$ (cf.\ Remark \ref{rem:scaling} for more details). 
The main result is that the distribution of $(1-\rho)Q$ converges to an exponential distribution, in line with classical heavy-traffic results that have been derived in a plethora of queueing models. 
As we will see, the parameter of the exponential distribution features the second moments of $\Lambda$ and $S$, which from now on we assume to exist. 

As a first result, however, we will show that, under a certain time-scaling, the marginal transient distributions of the process $(Q_n)_n$ converge to their reflected Brownian motion counterpart. 
We follow an argumentation developed in \cite[Chapter 5]{DM}.
We scale the queueing process by $(1-\rho)$ and time by a factor $(1-\rho)^{-2}$, in line with the usual heavy-traffic scaling. We wish to find the limit, as $\rho\uparrow 1$, for given $r,s>0$, of
\[K\left((1-\rho)^2r, e^{-(1-\rho)s}\right).\]
To this end, we first wish to identify $s_0(r, \rho)$ solving, in the regime $\rho\uparrow 1$,
\[e^{-(1-\rho)s} = \left(1-(1-\rho)^2r\right) \nu\left(e^{-(1-\rho)s}\right).\]
We do so by expanding both sides as polynomials in $1-\rho$, so as to obtain the following equation:
\begin{align*}1-(1-\rho)s +\tfrac{1}{2} (1-\rho)^2 s^2 &= \left(1-(1-\rho)^2r\right) \cdot\\
 \big(1-\nu'(1)\big((1-\rho)s& -\tfrac{1}{2}(1-\rho)^2 s^2\big)+\tfrac{1}{2}\nu''(1)(1-\rho)^2s^2\big)+O\big((1-\rho^3)\big).
\end{align*}
Using that $\rho=\nu'(1)$, we obtain the following equation after subtracting 1 from both sides and dividing by $(1-\rho)^2$:
\[\tfrac{1}{2}\nu''(1) \,s^2+s-r=O(1-\rho).\]
We thus conclude that
\[s_0(r) =\frac{-1+\sqrt{1+2\nu''(1)\,r}}{\nu''(1)} +O(1-\rho).\]
We assume that, as $\rho\uparrow 1$, $\kappa_0(e^{-(1-\rho)s})\to\bar\kappa_0(s)$ for some transform $\bar\kappa_0(\cdot)$;
this means that the initial distribution converges to some limiting distribution under the heavy-traffic scaling.
Now using the result stated in Theorem \ref{thtr}, we conclude after some standard computations 
that, as $\rho\uparrow 1$,
\[\ K\left((1-\rho)^2r, e^{-(1-\rho)s}\right) \to \frac{r}{r-s-\frac{1}{2}\nu''(1) s^2}\left(\bar \kappa_0(s)-\frac{s}{s_0(r)}\bar\kappa_0(s_0(r))\right).\]
In this expression we recognize the Laplace transform for the position of reflected Brownian motion after an exponentially distributed time (with mean $r^{-1}$), given the initial level has transform $\bar\kappa_0(\cdot)$; cf.\ for example \cite[Theorem 4.1]{DM}. A direct verification yields that 
$\nu''(1) =  {\mathbb E}N(N-1) = {\mathbb E}[\Lambda^2]\, {\mathbb E}[S^2]$.
Define $B(t)$ as $-t + \sqrt{\nu''(1)/2}\cdot W(t)$, with $W(\cdot)$ standard Brownian motion. In addition, let  $Q^{[x]}(\cdot)$ be the reflection of $B(\cdot)$ at 0, with initial condition $Q^{[x]}(0)=x.$
From the above argumentation, an application of the L\'evy convergence theorem provides us with the following result. 

\begin{theorem} \label{thhttr} Suppose $(1-\rho)Q_0= x\geqslant 0$.  Then we have for any $t>0$ that $(1-\rho)Q_{t/(1-\rho)^2}$ converges to $Q^{[x]}(t)$ as $\rho\uparrow 1$.
\end{theorem}

To study the stationary heavy-traffic behavior, we use a classical argumentation based on the Laplace-Stieltjes transform (\ref{LSTkappa}). Based on the fact that the stationary distribution of reflected Brownian motion is exponential, Theorem \ref{thhttr} suggest that $(1-\rho)Q$ converges in distribution to an exponentially distributed random variable with mean
 $\tfrac{1}{2} \ {\mathbb E}[\Lambda^2] \ {\mathbb E}[S^2]$. In the remainder of this subsection we make this claim precise. 

We evaluate $\kappa(e^{-(1-\rho)s})$, for $s\geqslant 0$ given, in the regime that $\rho\uparrow 1$. The starting point is
\[\nu(e^{-(1-\rho)s}) = \int_0^\infty \sigma\big(\lambda(1-e^{-(1-\rho)s})\big) \,g(\lambda)\,{\rm d}\lambda.\]
By expanding $e^{-(1-\rho)s}$, we first rewrite the expression from the previous display as
\[\int_0^\infty \sigma\left(\lambda(1-\rho)s-\frac{\lambda}{2}(1-\rho)^2s^2 +O\big((1-\rho^3) \big)\right) \,g(\lambda)\,{\rm d}\lambda.\]
Subsequently applying a Taylor expansion of $\sigma(\cdot)$ at 0, we obtain that $\nu(e^{-(1-\rho)s})$ can be expressed as
\[\int_0^\infty \left(1+\sigma'(0)\left(\lambda(1-\rho)s-\frac{\lambda}{2}(1-\rho)^2s^2\right) +
\frac{\sigma''(0)}{2}\lambda^2(1-\rho)^2 s^2+O\big((1-\rho^3) \big)
\right)\,g(\lambda)\,{\rm d}\lambda,\]
which can be interpreted as
\[1-(1-\rho)s\ {\mathbb E}\Lambda  \, {\mathbb E}S
+\tfrac{1}{2}(1-\rho)^2 s^2 \ {\mathbb E}\Lambda  \, {\mathbb E}S
 + \tfrac{1}{2}(1-\rho)^2 s^2\ {\mathbb E}[\Lambda^2] \ {\mathbb E}[S^2]+O\big((1-\rho^3) \big).\]
 To study the behavior in the heavy-traffic limit, we now consider the individual elements in the right-hand side of (\ref{LSTkappa}), evaluated in $e^{-(1-\rho)s}$. The numerator reads
 \begin{align*} 
 \nu(e^{-(1-\rho)s})\,\big(1-e^{-(1-\rho)s}\big) \left(1 -{{\mathbb E}\Lambda}\,{{\mathbb E}S}\right)&\sim (1-\rho)^2 s +O\big((1-\rho^3)\big),\end{align*}
 with `$\sim$' denoting that the ratio of the left- and right-hand side converges to 1 as $\rho\uparrow 1.$ The denominator expands as
 \[(1-\rho)^2s
 + \tfrac{1}{2}(1-\rho)^2 s^2\, {\mathbb E}[\Lambda^2]\,  {\mathbb E}[S^2]+O\big((1-\rho^3) \big).\]
 Upon combining the above, and by dividing both numerator and denominator by $(1-\rho)^2 s$, we obtain that, as $\rho\uparrow 1$,
 \[\kappa(e^{-(1-\rho)s}) \to \frac{1}{1 
 + \tfrac{1}{2} \ {\mathbb E}[\Lambda^2] \, {\mathbb E}[S^2]\,s}.\]
 We thus obtain the following result.
 \begin{theorem}\label{HT}As $\rho\uparrow 1$, we have that $(1-\rho)Q$ converges to an exponentially distributed random variable with mean
 $\tfrac{1}{2} \ {\mathbb E}[\Lambda^2] \, {\mathbb E}[S^2]$.
 \end{theorem}
 
 \begin{remark}{\em Theorem \ref{HT} covers the heavy-traffic distribution in the ordinary M/G/1 queue as a special case. In that case the arrival rate is deterministic, such that ${\mathbb E}[\Lambda^2]= ({\mathbb E}\Lambda)^2$. With a bit of rewriting, we find that in this case $(1-\rho)Q$ converges to an exponentially distributed random variable with mean $\frac12 \E[S^2]/(\E S)^2$.}
\end{remark}

\section{Conclusions and suggestions for further research} \label{sec:conclusion} 
In this paper we have considered various single-server queues under overdispersion. 
While their infinite-server counterparts allow for explicit analysis, these queues can rarely be dealt with in a straightforward manner. 
Focusing on M/M/1 queues with resampled arrival and service rates, closed-form expressions for the queue length can only be derived if these rates can attain just finitely many values (Section \ref{sec:exact}). 
In the heavy-traffic regime, however, explicit convergence results can be derived; specifically, we have shown that a scaled version of the steady-state queue length converges to an exponentially distributed random variable (Sections \ref{sec:exo} and \ref{sec:exoperiodic}). 
We also considered a model in which the arrival rate is resampled upon service completion (rather than at i.i.d.\ resampling times, independently of the queue's dynamics); also in this setting heavy-traffic analysis has been performed (Section \ref{sec:endo}).

In the area of queues under overdispersion there are still many open problems. 
The most prominent question is posed in Section \ref{sec:exact}: can we find the queue-length distribution if the arrival and service rates are repeatedly sampled from a distribution with countably infinite or even uncountable support? 
As we demonstrated, the conventional queueing-theoretic approach, expressing the corresponding Laplace transform  through a system of finitely many equations with equally many unknowns, clearly does not apply.

As overdispersion was observed in various types of service systems  \cite{BAS, LIU,WHI2} that are typically modeled as {\it many-server} queues (such as call centers),
one would like to get a handle on such queues with resampled rates as well. 
The ultimate goal would be to design staffing rules for many-server queues under overdispersion; see e.g. \cite{HLMM20}. 

\subsection*{Acknowledgments}
The authors would like to thank Bo Klaasse and Rudesindo N\'u\~nez-Queija for useful comments and suggestions.


\begin{appendices}
\section{Appendix: alternative computation of the asymptotic variance}\label{APP}
The asymptotic (co-)variances $v_A$, $v_S$, and $c_{A,S}$ can be computed in an alternative way, using results from large deviations theory \cite{DOC}; a similar approach has been followed in e.g.\ \cite{SBM,KM}. With this approach, also higher (centered) moments of $A(t)/t$ and $S(t)/t$ can be calculated in closed form in the regime that $t\to\infty$, as we point out below.

Let $X(\cdot)$ be a stochastic process with stationary increments, and let $c$ be larger than 
$\mathbb{E} X(1)$. Define the {\it asymptotic cumulant generating function}
\[
\gamma(\theta) := \lim_{t \to \infty} \frac{1}{t} \log \mathbb{E} \,\mathrm{exp}(\theta X(t)),
\]
assumed to be finite in an open neighborhood of the origin.
Then,  according to \cite{DOC}, as $u \to \infty$,
\[
\frac{1}{u} \log \mathbb{P}(\exists t>0: X(t) - ct \geq u) \to -\theta^\star,
\]
where $\theta^\star$ solves $c \theta = \gamma(\theta)$, or
\[ \lim_{t \to \infty} \frac{1}{t} \log \mathbb{E} \,\mathrm{exp}(\theta X(t)) - c\theta=0.\]

Interestingly, in case there is a regenerative structure underlying the process $X(\cdot)$, there is a second way of computing $\theta^\star.$
Let $T_n$ be time epochs such that $X(T_n)-X(T_{n-1})$ are i.i.d., and let $T:=T_1$. Then  we equivalently have
\[
\frac{1}{u} \log \mathbb{P}(\exists n>0: X(T_n) - cT_n \geq u) \to -\theta^\star, 
\]
where (the same!) $\theta^\star$ solves
\begin{equation}\label{DufCon}
\log\mathbb{E} \exp(\theta (X(T) - c T)) = 0.
\end{equation}
The idea is that we are going to exploit the fact that, obviously, both procedures should lead to  the same $(c,\theta)$ pairs.

In our setup, we identify $\xi$ with $T$: the regeneration points are the resample epochs. 
For reasons that will become clear later, 
we study the process $X_{\alpha}(t)$ that is a linear combination of the arrival and service process: 
we define $X_{\alpha}(t):=\alpha A(t) + (1-\alpha) S(t)$ for $\alpha\in{\mathbb R}$.
We show how to compute
\[v_X^{(\alpha)} =\lim_{t\to\infty}\frac{{\mathbb V}{\rm ar}\,X_{\alpha}(t)}{t},\]
which can be checked to equal $\gamma''(0).$ 
In the sequel, we study the function $c(\theta) := \gamma(\theta)/\theta$; as both $c(\cdot)$ and $\gamma(\cdot)$ depend on $\alpha$, we will consistently write $c_\alpha(\cdot)$ and $\gamma_\alpha(\cdot)$.
Once having found $v_X^{(\alpha)}$, it turns out that with the right choices of $\alpha$ we can compute $v_A$, $v_S$, and $c_{A,S}$.

Let $\Xi(\cdot)$ be the moment generating function of the resampling time $\xi$. It takes an elementary calculation (by conditioning on the value of $\Lambda$ and ${\mathcal M}$ in the interval under consideration, as well as on the duration of the inter-sample interval) to verify that
\begin{align*}
\varphi_\alpha(c,\theta):=\mathbb{E} \exp({\theta (X_\alpha (T)  - c T)} )&=\mathbb{E} \exp({\theta (\alpha A(T) + (1-\alpha) S(T) - c T)} )\\
&= \mathbb{E}\left[\Xi\big( \Lambda(\mathrm{e}^{\theta \alpha}-1) + \mathcal{M}(\mathrm{e}^{\theta (1-\alpha)}-1) - c \theta\big)\right].
\end{align*}
Applying $(\ref{DufCon})$, we find that $c_\alpha(\theta)$ is the value of $c$ for which $\log \varphi_\alpha(c,\theta) =0;$
differentiating the equivalent relation $\varphi_\alpha(c_\alpha(\theta),\theta) =1$ with respect to $\theta$ gives
\begin{align*}
\frac{\mathrm{d}}{\mathrm{d}\theta} \varphi_\alpha(c_\alpha(\theta),\theta)
&= \frac{\mathrm{d}}{\mathrm{d}\theta}\mathbb{E}\left[ \Xi\big( \Lambda(\mathrm{e}^{\theta \alpha}-1) + \mathcal{M}(\mathrm{e}^{\theta (1-\alpha)}-1) - c_\alpha(\theta) \theta\big)\right]\\
&= \mathbb{E}\, \Xi' (\Lambda(\mathrm{e}^{\theta \alpha}-1) + \mathcal{M}(\mathrm{e}^{\theta (1-\alpha)}-1) - c_\alpha(\theta) \theta) \\
& \qquad \qquad \cdot \left[\alpha \Lambda \mathrm{e}^{\theta \alpha} + (1-\alpha) \mathcal{M}\mathrm{e}^{\theta (1-\alpha)} - c_\alpha(\theta) - \theta c_\alpha'(\theta) \right]
= 0.
\end{align*}
Inserting $\theta=0$ gives the (obvious, but reassuring) identity
\begin{equation}\label{cknot}
\lim_{t\to\infty}\frac{ {\mathbb E}X_\alpha(t)}{t} =c_\alpha(0)=\alpha\, \E \Lambda + (1-\alpha) \,\E \mathcal{M}.
\end{equation}
Again differentiating the above equation gives
\begin{align*}
\frac{\mathrm{d^2}}{\mathrm{d}\theta^2} \varphi_\alpha(c_\alpha(\theta),\theta)=&\: \mathbb{E}\, \Xi''(\Lambda(\mathrm{e}^{\theta \alpha}-1) + \mathcal{M}(\mathrm{e}^{\theta (1-\alpha)}-1) - c_\alpha(\theta) \theta)\\
&\qquad  \cdot \left[\alpha \Lambda \mathrm{e}^{\theta \alpha} + (1-\alpha) \mathcal{M}\mathrm{e}^{\theta (1-\alpha)} - c_\alpha(\theta) - \theta c_\alpha'(\theta) \right]^2\:+\\
&  \mathbb{E} \,\Xi'(\Lambda(\mathrm{e}^{\theta \alpha}-1) + \mathcal{M}(\mathrm{e}^{\theta (1-\alpha)}-1) - c_\alpha(\theta) \theta)\\
&\qquad \cdot \left[\alpha^2 \Lambda \mathrm{e}^{\theta \alpha} + (1-\alpha)^2 \mathcal{M}\mathrm{e}^{\theta (1-\alpha)} - 2c_\alpha'(\theta) - \theta c_\alpha''(\theta) \right]= 0.
\end{align*}
Evaluating this equality in $\theta=0$ gives 
\[\E \xi^2 \cdot \E[\alpha  \Lambda + (1-\alpha)\mathcal{M}-c_\alpha(0)]^2 + \E \xi \cdot(\alpha^2 \E\Lambda + (1-\alpha)^2 \E\mathcal{M} - 2c_\alpha'(0))=0,\]
or, equivalently, using $(\ref{cknot})$,
\begin{equation} \label{BAS}
v_X^{(\alpha)}=\gamma_\alpha''(0)=2c_\alpha'(0) =\frac{ \E \xi^2}{\E\xi} \cdot {\mathbb V}{\rm ar}(\alpha  \Lambda + (1-\alpha)\mathcal{M}) +  (\alpha^2 \E\Lambda + (1-\alpha)^2 \E\mathcal{M}).
\end{equation}
We are now in a position to compute $v_A$, $v_S$, and $c_{A,S}$. 
To this end, first realize that inserting $\alpha=0$, $\alpha=1$ and $\alpha=\tfrac{1}{2}$ yields, respectively,
\[
v_{X}^{(1)} = v_A,\:\:\:\:v_{X}^{(0)} = v_S,\:\:\:\:v_{X}^{(1/2)} = \frac14 v_A + \frac14 v_S + \frac12 c_{A,S},
\]
where, by virtue of (\ref{BAS}), 
\begin{align*}
v_{X}^{(1)}&= {\mathbb E}\Lambda+{{\mathbb V}{\rm ar}\,\Lambda}\,\frac{{\mathbb E}[\xi^2]}{{\mathbb E}\xi},\:\:\:\:v_{X}^{(0)}= {\mathbb E}\mathcal M+{{\mathbb V}{\rm ar}\,\mathcal M}\,\frac{{\mathbb E}[\xi^2]}{{\mathbb E}\xi},\\
v_{X}^{(1/2)}&=  \tfrac14 \E \Lambda + \tfrac14 \E \mathcal{M} + \frac{{\mathbb E}[\xi^2]}{{\mathbb E}\xi} {\mathbb V}{\rm ar} (\tfrac12 \Lambda + \tfrac12 \mathcal{M})\frac{{\mathbb E}[\xi^2]}{{\mathbb E}\xi}.
\end{align*}
The claim of  Corollary \ref{C1} now follows immediately.

The technique presented here has a clear advantage over the direct approach that was used at the end of Section $\ref{arpos}$:
it can be used to find higher cumulant moments. 
As is readily checked, for instance for the third asymptotic cumulant moment
\[\lim_{t\to\infty}\frac{\E(X_\alpha(t)- \E X_\alpha(t))^3}{t}=\gamma_\alpha'''(0)=3c_\alpha''(0).\]
The value of $c_\alpha''(0)$ can be found by differentiating the relation $\varphi_\alpha(c_\alpha(\theta),\theta) =1$ three times with respect to $\theta$; the corresponding calculations are straightforward yet tedious.

\end{appendices}


\bibliographystyle{plain}
\end{document}